\documentclass{amsart}
\usepackage{amsmath,amssymb,latexsym,enumerate}
\usepackage{pstricks}
\usepackage{pst-node}
\usepackage{tikz}
\usepackage{imakeidx}

\newtheorem{theorem}{Theorem}[section]

\newtheorem{lemma}[theorem]{Lemma}
\newtheorem{proposition}[theorem]{Proposition}

\theoremstyle{definition}
\newtheorem{definition}[theorem]{Definition}
\newtheorem{example}[theorem]{Example}

\newtheorem{subexercise}{Exercise}[theorem]
\newtheorem{exercise}[theorem]{Exercise}
\theoremstyle{remark} 
\newtheorem{remark}[theorem]{Remark}
\numberwithin{equation}{section}

\newcommand{\cA}{{\mathcal{A}}}
\newcommand{\cF}{{\mathcal{F}}}
\newcommand{\cL}{{\mathcal{L}}}
\newcommand{\cO}{{\mathcal{O}}}
\newcommand{\cT}{{\mathcal{T}}}
\newcommand{\cU}{{\mathcal{U}}}
\newcommand{\ba}{{\bf a}}
\newcommand{\bb}{{\bf b}}

\newcommand{\x}{{\bf x}}
\newcommand{\bx}{{\bf x}}
\newcommand{\y}{{\bf y}}

\newcommand{\CC}{{\mathbb{C}}}
\newcommand{\PP}{{\mathbb{P}}}
\newcommand{\ZZ}{{\mathbb{Z}}}

\newcommand{\meet}[2]{#1 \cap #2}
\newcommand{\join}[2]{\langle #1, #2 \rangle}
\newcommand{\sgn}{\operatorname{sign}}

\newcommand{\erase}[1]{{}}

\makeindex

\title{Introduction to Cluster Algebras}
\author{Max Glick}
\author{Dylan Rupel}

\begin{document}
\begin{abstract}
  These are notes for a series of lectures presented at the ASIDE conference 2016.  The definition of a cluster algebra is motivated through several examples, namely Markov triples, the Grassmannians $Gr_2(\CC^n)$, and the appearance of double Bruhat cells in the theory of total positivity.  Once the definition of cluster algebras is introduced in several stages of increasing generality, proofs of fundamental results are sketched in the rank 2 case.  From these foundations we build up the notion of Poisson structures compatible with a cluster algebra structure and indicate how this leads to a quantization of cluster algebras.  Finally we give applications of these ideas to integrable systems in the form of Zamolodchikov periodicity and the pentagram map.
\end{abstract}
\maketitle

\section{Introduction}
  Cluster algebras were introduced by Fomin and Zelevinsky \cite{FZ02} in 2002 as the culmination of their study of total positivity \cite{FZ99} and (dual) canonical bases.  The topic of cluster algebras quickly grew into its own as a subject deserving independent study mainly fueled by its emergent close relationship to many areas of mathematics.  Here is a partial list of related topics: combinatorics \cite{MP07}, hyperbolic geometry \cite{FG06,FST08,MSW13}, Lie theory \cite{GLS06}, Poisson geometry \cite{GSV03}, integrable systems \cite{dFK10,G11}, representations of associative algebras \cite{CC06,CK06,BMRRT06,R11,Q12,R15}, mathematical physics \cite{EF12,ABCGPT14}, and quantum groups \cite{K12,GLS13,KQ14,BR15}.

  In these notes we will give an introduction to cluster algebras and a couple of the applications mentioned above.  These notes are far from exhaustive and the references above only touch on the vast literature.  Other overviews of cluster algebras can be found in the works \cite{K10,W13,GLS} which will also provide additional references.

  These lecture notes are organized as follows.  Section~\ref{sec:motivation} gives several motivating examples from which we will abstract the definition of a cluster algebra.  
	Section~\ref{sec:cluster_algebras} contains several variations on the definition of cluster algebras with increasing generality and also the definition of $Y$-patterns.  In Section~\ref{sec:basic_theorems}, we describe the foundational results in the theory of cluster algebras and sketch, or otherwise indicate the ideas behind, the proofs of these results. Section~\ref{sec:poisson_and_quantum} recalls the theory of Poisson structures compatible with a cluster algebra and describes how this naturally leads to a quantization of cluster algebras.  Finally, we conclude with applications of the cluster algebra machinery to problems involving integrable systems in Section~\ref{sec:integrable_systems}.

\section{Motivating Examples}\label{sec:motivation}
	Cluster algebras are certain commutative rings possessing additional structure, including:
	\begin{itemize}
		\item a distinguished collection of generators called \emph{cluster variables};
		\item a collection of finite subsets of the set of all cluster variables called \emph{clusters};
		\item a \emph{mutation} rule which, given a cluster and one of its variables $x$, produces another cluster by replacing $x$ with a different cluster variable $x'$ related to $x$ by an \emph{exchange relation}
		\begin{displaymath}
			xx' = F,
		\end{displaymath}
		where $F$ is a binomial in the variables common to the initial and mutated clusters.
	\end{itemize}
	Before defining cluster algebras, we describe a few settings illustrating these various components.  We leave it to the astute reader to compare these examples with the definitions presented in Section~\ref{sec:cluster_algebras} and to see how they may be reinterpreted within the cluster algebra framework.

	\subsection{Markov Triples} \label{sub:Markov triples}
    A \emph{Markov triple}\index{Markov!1@triple} is a triple $(a,b,c)$ of positive integers satisfying the \emph{Markov equation} $a^2+b^2+c^2=3abc$; an integer which appears as a term in a Markov triple is called a \emph{Markov number}.  The Markov equation is an example of a Diophantine equation and two classical number theoretic problems are to determine the number of solutions and to determine a method for finding all such solutions.  We will solve both of these problems for the Markov equation.
    
    \begin{exercise}
      Prove that $(1,1,1)$ and $(1,1,2)$ are (up to reordering) the only Markov triples with repeated values.
    \end{exercise}
    Rearranging the Markov equation we see that $c^2-3abc+a^2+b^2=0$ and so $c$ is a root of the quadratic $f(x)=x^2-(3ab)x+(a^2+b^2)$.  But notice that the other root $c'=\frac{a^2+b^2}{c}=3ab-c$ is a positive integer and so $(a,b,c')$ is again a Markov triple.  

    Note that there was nothing special about $c$ in the calculation above.  Thus given any Markov triple $(a,b,c)$, we may perform three possible exchanges 
    \[\Big(a,b,\frac{a^2+b^2}{c}\Big)\quad\Big(a,\frac{a^2+c^2}{b},c\Big)\quad\Big(\frac{b^2+c^2}{a},b,c\Big)\]
    and obtain another Markov triple in each case.  The following exercises solve the above two classical problems of Diophantine equations.
    \begin{exercise}\mbox{}
      \begin{enumerate}[\quad\upshape (a)]
        \item Prove that there are infinitely many Markov triples by showing that there is no bound on how large the largest value can be.
        \item Show that all Markov triples may be obtained from the Markov triple $(1,1,1)$ by a sequence of exchanges.
      \end{enumerate}
    \end{exercise}

	\subsection{The Grassmannian $Gr_{2}(\CC^n)$} \label{sub:Gr2n}
		The \emph{Grassmannian}\index{Grassmannian} $Gr_{k}(\CC^n)$ is the set of $k$-dimensional linear subspaces of $\CC^n$.  A point in the Grassmannian can be described, albeit non-uniquely, as the row span of a full rank matrix $A \in \CC^{k \times n}$.  The maximal minors of $A$ are called \emph{Pl\"ucker coordinates}\index{Pl\"ucker!coordinate}.
		
		Now, restrict to $k=2$ and let 
		\begin{displaymath}
			A = \left[ \begin{array}{cccc}
			a_{11} & a_{12} & \ldots & a_{1n} \\
			a_{21} & a_{22} & \ldots & a_{2n} \\
			\end{array} \right].
		\end{displaymath}
		The Pl\"ucker coordinates are given by $\Delta_{ij} = a_{1i}a_{2j} - a_{1j}a_{2i}$ for $1 \leq i < j \leq n$.  As the matrix $A$ is not uniquely determined by a point in $Gr_2(\CC^n)$, the Pl\"ucker coordinates are not truly functions on $Gr_2(\CC^n)$, but rather they are only well-defined functions up to simultaneous rescaling.  That is, the Pl\"ucker coordinates determine a closed embedding into the projective space $\PP^{{n\choose 2}-1}$.  The following exercise shows that the image of the Pl\"ucker embedding lies inside a closed subset of projective space.
		\begin{exercise}
			Verify that the $\Delta_{ij}$ satisfy the so-called \emph{Pl\"ucker relations}\index{Pl\"ucker!relation}
			\begin{equation}\label{eq:short plucker special}
				\Delta_{ik}\Delta_{jl} = \Delta_{ij}\Delta_{kl} + \Delta_{il}\Delta_{jk}
			\end{equation}
			for $1 \leq i < j < k < l \leq n$.
		\end{exercise}
		
		Consider a regular $n$-gon with vertices labeled $1,2,\ldots, n$.  Let $T$ be a triangulation\index{triangulation}, i.e. a maximal collection of chords $\overline{ij}$ with $1 \leq i < j \leq n$, no two of which intersect in their interiors.  Note that $T$ always consists of $n-3$ diagonals together with the $n$ sides $\overline{12}, \overline{23}, \overline{34}, \ldots, \overline{1n}$.  Associate to $T$ the corresponding collection of Pl\"ucker coordinates
		\begin{displaymath}
			\Delta_T := \{\Delta_{ij} : \overline{ij} \in T\}.
		\end{displaymath}
    \setcounter{subproposition}{1}
		\begin{proposition} \label{prop:Gr2n}
		Fix positive reals $x_{ij}$ for all $\overline{ij} \in T$.  Then there exists $A \in Gr_2(\CC^n)$ such that $\Delta_{ij}(A) = x_{ij}$ for all $\overline{ij} \in T$.  Moreover, each $\Delta_{ij}(A)$ with $\overline{ij} \notin T$ can be expressed as a subtraction free rational expression of the $x_{ij}$.
		\end{proposition}
		
		For instance, let $n=4$ and $T = \{\overline{13},\overline{12},\overline{23},\overline{34},\overline{14}\}$.  A possible representing matrix is
		\begin{displaymath}
			A = \left[\begin{array}{cccc}
			1 & \frac{x_{23}}{x_{13}} & 0 & -\frac{x_{34}}{x_{13}} \\
			0 & x_{12} & x_{13} & x_{14} \\
			\end{array}\right].
		\end{displaymath}
		There is one remaining Pl\"ucker coordinate, namely $\Delta_{24}$, and it can be computed from the given ones as
		\begin{displaymath}
		\Delta_{24}(A) = \frac{x_{12}x_{34} + x_{14}x_{23}}{x_{13}}
		\end{displaymath}
		which follows from a Pl\"ucker relation.	This formula can be thought of as a change of coordinates from $\Delta_T$ to $\Delta_{T'}$ where $T' = (T \setminus \{\overline{13}\}) \cup \{\overline{24}\}$.  
    \setcounter{subexercise}{2}
    \begin{exercise}
      The change of coordinate systems observed above can be understood more classically in terms of the Ptolemy relations for cyclic quadrilaterals.  Indeed, consider four distinct points labeled $1,2,3,4$ on a circle in cyclic order and let $x_{ij}$ denote the distance between vertices $i$ and $j$.  Prove that $x_{13}x_{24}=x_{12}x_{34}+x_{14}x_{23}$.
    \end{exercise}
		
		More generally, if $T$ is a triangulation of an $n$-gon and $\overline{ik} \in T$ is a diagonal (as opposed to a side), then $\overline{ik}$ is part of two triangles of $T$.  Call the third vertices of these two triangles $j$ and $l$.  It follows that $T' = (T \setminus \{\overline{ik}\}) \cup \{\overline{jl}\}$ is again a triangulation.  Moreover, it is well known that any two triangulations may be related by a sequence of these flips of diagonals \cite{H91}.  Performing a sequence of such quadrilateral flips and using a Pl\"ucker relation at each step makes it possible to iteratively compute all the rational expressions promised in Proposition \ref{prop:Gr2n}.
 
	\subsection{Double Bruhat cells} \label{sub:total positivity}
    An $n\times n$ matrix $M$ is called \emph{totally positive}\index{totally positive matrix} if the determinant of every square submatrix is a positive real number.  In particular, every entry of $M$ is positive and $M$ must be invertible.  Write $GL_n^{>0}\subset GL_n$ for the subset of totally positive matrices.  To check that a given matrix $M\in GL_n$ is totally positive one must, a priori, check that all ${2n\choose n}-1$ minors of $M$ are positive.  A natural question is whether this verification process can be made more efficient.  More precisely, is there a smaller collection of minors one may compute and from the positivity of this subset conclude that every minor is positive, i.e. conclude that $M\in GL_n^{>0}$?  We will call such a collection a \emph{total positivity criterion}\index{total positivity criterion} if it exists.  

    For small $n$, such criteria can be found and verified easily.  For example, a matrix $M=\left[\begin{array}{cc} a & b\\ c & d\end{array}\right]\in GL_2$ is totally positive if and only if $a,b,c,ad-bc>0$ if and only if $b,c,d,ad-bc>0$.
    \begin{exercise}\label{exer:gl3 positivity}
      Find a minimal collection of minors whose positivity guarantees a matrix in $GL_3$ is totally positive (hint: any such total positivity criterion consists of 9 minors).
    \end{exercise}

    To describe a solution and easily identify total positivity criteria for all general linear groups $GL_n$, it will be convenient to slightly generalize the notion of total positivity.  An $n\times n$ matrix $M$ is called \emph{totally nonnegative}\index{totally nonnegative matrix} if the determinant of every square submatrix is a nonnegative real number.  Write $GL_n^{\ge0}\subset GL_n$ for the subset of totally nonnegative matrices.  Again one may ask: what is a minimal collection of minors needed to check that a matrix is totally nonnegative?  Unfortunately, or perhaps fortunately, the total nonnegativity criteria are not uniformly described across all of $GL_n$.  The solution to this problem naturally leads one to study certain subvarieties of $GL_n$ called \emph{double Bruhat cells}, which we now describe.

    Let $B_+,B_-\subset GL_n$ denote the subgroups of upper and lower triangular matrices respectively.  Identify the symmetric group $\Sigma_n$ with the subgroup of $GL_n$ consisting of permutation matrices, i.e. matrices having precisely one nonzero entry 1 in each row and column.  For example, identify the permutation $(1\ 2)\in\Sigma_2$ (written in cycle notation) with the matrix $\left[\begin{array}{cc}0 & 1\\ 1 & 0\end{array}\right]$.  

    It is well known that $GL_n$ decomposes in two ways (actually many ways) as a union of \emph{Bruhat cells}:
    \[GL_n=\bigsqcup_{w\in\Sigma_n}B_+w B_+=\bigsqcup_{w\in\Sigma_n}B_-w B_-.\]
    To understand total nonnegativity criteria for $GL_n$, we will consider the \emph{double Bruhat cells}\index{double Bruhat cell} $GL_n^{u,v}=B_+uB_+\cap B_-vB_-$ for $u,v\in\Sigma_n$.  
    \setcounter{subexercise}{2}
    \begin{exercise}
      Find all 4 double Bruhat cells in $GL_2$ and verify that they partition the space.  If you are brave, find all double Bruhat cells in $GL_3$ and verify that they partition the space.
    \end{exercise}
    For $I,J\subset[1,n]$ with $|I|=|J|$, denote by $\Delta_{I,J}$ the function on $GL_n$ which returns the determinant of the submatrix on row set $I$ and column set $J$.
    \setcounter{subexample}{1}
    \begin{example}\label{example:big cell}
      Let $w_0=(1\ n)(2\ n-1)\cdots(\lfloor n/2\rfloor\ \lceil n/2\rceil)$ denote the longest permutation in $\Sigma_n$.  The double Bruhat cell $GL_n^{w_0,w_0}$ is given as follows:
      \[GL_n^{w_0,w_0}=\{M\in GL_n:\Delta_{[1,i],[n+1-i,n]}(M)\ne0,\Delta_{[n+1-i,n],[1,i]}(M)\ne0\text{ for all $i$}\}.\]
    \end{example}

\begin{figure}[ht]
    \begin{tikzpicture}
      \draw (0,0) -- (0.5,0) -- (1.5,1) -- (2,1) -- (3,2) -- (3.5,2) -- (4.5,3) -- (5,3) -- (6,3) -- (6.5,3) -- (7.5,3) -- (8,3) -- (9,3) -- (9.5,3);
      \draw (0,1) -- (0.5,1) -- (1.5,0) -- (2,0) -- (3,0) -- (3.5,0) -- (4.5,0) -- (5,0) -- (6,1) -- (6.5,1) -- (7.5,2) -- (8,2) -- (9,2) -- (9.5,2);
      \draw (0,2) -- (0.5,2) -- (1.5,2) -- (2,2) -- (3,1) -- (3.5,1) -- (4.5,1) -- (5,1) -- (6,0) -- (6.5,0) -- (7.5,0) -- (8,0) -- (9,1) -- (9.5,1);
      \draw (0,3) -- (0.5,3) -- (1.5,3) -- (2,3) -- (3,3) -- (3.5,3) -- (4.5,2) -- (5,2) -- (6,2) -- (6.5,2) -- (7.5,1) -- (8,1) -- (9,0) -- (9.5,0);
    \end{tikzpicture}
    \caption{A wiring diagram for the reduced word $(1,2,3,1,2,1)\in\Sigma_4$.}
    \label{fig:wiring diagram}
    \end{figure}
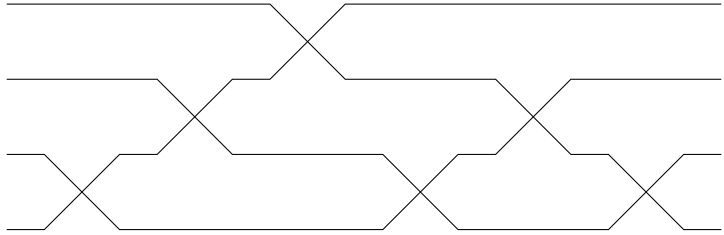

    It turns out that each double Bruhat cell admits its own collection of total positivity criteria, i.e. for each $u,v\in\Sigma_n$ there exists a minimal collection of minors whose positivity identifies the subset $GL_n^{u,v}\cap GL_n^{\ge0}$ inside the double Bruhat cell $GL_n^{u,v}$.
    Some of these total nonnegativity criteria in $GL_n^{u,v}$ can be conveniently described using double wiring diagrams.  Write $s_i=(i\ i+1)\in\Sigma_n$ for the simple transposition interchanging $i$ and $i+1$.  A \emph{reduced word}\index{reduced word} for $w\in\Sigma_n$ is a minimal sequence $(i_1,\ldots,i_r)$ so that $w=s_{i_1}\cdots s_{i_r}$, where $\ell(w):=r$ is called the \emph{length} of $w$.  A reduced word $(i_1,\ldots,i_r)$ is naturally encoded in a \emph{wiring diagram}\index{wiring diagram}, i.e. a collection of $n$ strands with a crossing between the $i$th and $(i+1)$th strands each time $i$ appears in the reduced word, see Figure~\ref{fig:wiring diagram}.

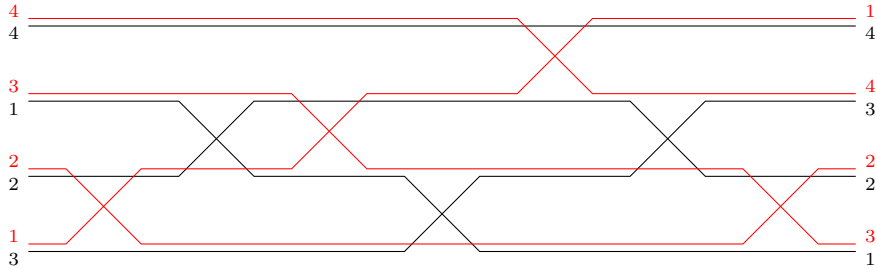
    \begin{figure}[ht]
    \begin{tikzpicture}
      \draw (0,0) -- (2,0) -- (3,0) -- (5,0) -- (6,1) -- (8,1) -- (9,2) -- (11,2);
      \draw (0,1) -- (2,1) -- (3,2) -- (5,2) -- (6,2) -- (8,2) -- (9,1) -- (11,1);
      \draw (0,2) -- (2,2) -- (3,1) -- (5,1) -- (6,0) -- (8,0) -- (9,0) -- (11,0);
      \draw (0,3) -- (2,3) -- (3,3) -- (5,3) -- (6,3) -- (8,3) -- (9,3) -- (11,3);
      \draw[color=red] (0,0.1) -- (0.5,0.1) -- (1.5,1.1) -- (3.5,1.1) -- (4.5,2.1) -- (6.5,2.1) -- (7.5,3.1) -- (9.5,3.1) -- (10.5,3.1) -- (11,3.1);
      \draw[color=red] (0,1.1) -- (0.5,1.1) -- (1.5,0.1) -- (3.5,0.1) -- (4.5,0.1) -- (6.5,0.1) -- (7.5,0.1) -- (9.5,0.1) -- (10.5,1.1) -- (11,1.1);
      \draw[color=red] (0,2.1) -- (0.5,2.1) -- (1.5,2.1) -- (3.5,2.1) -- (4.5,1.1) -- (6.5,1.1) -- (7.5,1.1) -- (9.5,1.1) -- (10.5,0.1) -- (11,0.1);
      \draw[color=red] (0,3.1) -- (0.5,3.1) -- (1.5,3.1) -- (3.5,3.1) -- (4.5,3.1) -- (6.5,3.1) -- (7.5,2.1) -- (9.5,2.1) -- (10.5,2.1) -- (11,2.1);
      \node[left] at (0,-0.1) {$\scriptstyle 3$};
      \node[left] at (0,0.9) {$\scriptstyle 2$};
      \node[left] at (0,1.9) {$\scriptstyle 1$};
      \node[left] at (0,2.9) {$\scriptstyle 4$};
      \node[left, color=red] at (0,0.2) {$\scriptstyle 1$};
      \node[left, color=red] at (0,1.2) {$\scriptstyle 2$};
      \node[left, color=red] at (0,2.2) {$\scriptstyle 3$};
      \node[left, color=red] at (0,3.2) {$\scriptstyle 4$};
      \node[right] at (11,-0.1) {$\scriptstyle 1$};
      \node[right] at (11,0.9) {$\scriptstyle 2$};
      \node[right] at (11,1.9) {$\scriptstyle 3$};
      \node[right] at (11,2.9) {$\scriptstyle 4$};
      \node[right, color=red] at (11,0.2) {$\scriptstyle 3$};
      \node[right, color=red] at (11,1.2) {$\scriptstyle 2$};
      \node[right, color=red] at (11,2.2) {$\scriptstyle 4$};
      \node[right, color=red] at (11,3.2) {$\scriptstyle 1$};
    \end{tikzpicture}\\
    \[\Delta_{1,3},\Delta_{2,3},\Delta_{2,1},\Delta_{3,1},\Delta_{12,23},\Delta_{12,13},\Delta_{23,13},\Delta_{23,12},\Delta_{123,123},\Delta_{234,123},\Delta_{1234,1234}\]
    \caption{A double wiring diagram for the double reduced word $(-1,2,-2,1,-3,2,-1)\in\Sigma_4$ and the collection of chamber minors determined by this double wiring diagram.}
    \label{fig:double wiring diagram}
    \end{figure}

    Given $(u,v)\in\Sigma_n^2$, a double reduced word for $(u,v)$ is an arbitrary shuffle of a reduced word for $u$ and a reduced word for $v$, where terms in the reduced word for $v$ are taken from the set $\{-1,\ldots,-(n-1)\}$ for clarity.  For example, taking $u=(1\ 3\ 2),v=(1\ 3)\in\Sigma_3$, one double reduced word for $(u,v)$ is $(-1,2,-2,1,-1)$.  Given a double reduced word $(i_1,\ldots,i_r)$, build a double wiring diagram\index{double wiring diagram} by superposing wiring diagrams associated to the reduced word for $u$ and the reduced word for $v$ with crossings ordered according to the double reduced word (see Figure~\ref{fig:double wiring diagram}).  

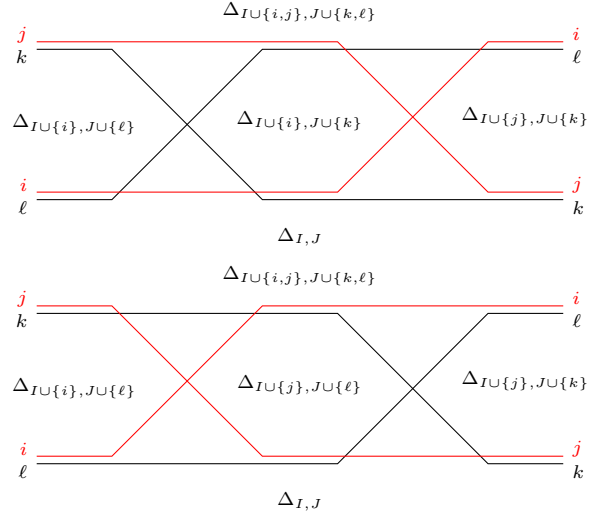
\begin{figure}[ht]
      \begin{tikzpicture}
        \draw (0,0) -- (1,0) -- (3,2) -- (7,2);
        \draw (0,2) -- (1,2) -- (3,0) -- (7,0);
        \draw[color=red] (0,0.1) -- (4,0.1) -- (6,2.1) -- (7,2.1);
        \draw[color=red] (0,2.1) -- (4,2.1) -- (6,0.1) -- (7,0.1);
        \node[left] at (0,-0.1) {$\scriptstyle \ell$};
        \node[left] at (0,1.9) {$\scriptstyle k$};
        \node[left, color=red] at (0,0.2) {$\scriptstyle i$};
        \node[left, color=red] at (0,2.2) {$\scriptstyle j$};
        \node[right] at (7,-0.1) {$\scriptstyle k$};
        \node[right] at (7,1.9) {$\scriptstyle \ell$};
        \node[right, color=red] at (7,0.2) {$\scriptstyle j$};
        \node[right, color=red] at (7,2.2) {$\scriptstyle i$};
        \node at (3.5,2.5) {$\scriptstyle\Delta_{I\cup\{i,j\},J\cup\{k,\ell\}}$};
        \node at (6.5,1.1) {$\scriptstyle\Delta_{I\cup\{j\},J\cup\{k\}}$};
        \node at (3.5,1.05) {$\scriptstyle\Delta_{I\cup\{i\},J\cup\{k\}}$};
        \node at (0.5,1) {$\scriptstyle\Delta_{I\cup\{i\},J\cup\{\ell\}}$};
        \node at (3.5,-0.5) {$\scriptstyle\Delta_{I,J}$};
      \end{tikzpicture}\\
      \begin{tikzpicture}
        \draw (0,0) -- (4,0) -- (6,2) -- (7,2);
        \draw (0,2) -- (4,2) -- (6,0) -- (7,0);
        \draw[color=red] (0,0.1) -- (1,0.1) -- (3,2.1) -- (7,2.1);
        \draw[color=red] (0,2.1) -- (1,2.1) -- (3,0.1) -- (7,0.1);
        \node[left] at (0,-0.1) {$\scriptstyle \ell$};
        \node[left] at (0,1.9) {$\scriptstyle k$};
        \node[left, color=red] at (0,0.2) {$\scriptstyle i$};
        \node[left, color=red] at (0,2.2) {$\scriptstyle j$};
        \node[right] at (7,-0.1) {$\scriptstyle k$};
        \node[right] at (7,1.9) {$\scriptstyle \ell$};
        \node[right, color=red] at (7,0.2) {$\scriptstyle j$};
        \node[right, color=red] at (7,2.2) {$\scriptstyle i$};
        \node at (3.5,2.5) {$\scriptstyle\Delta_{I\cup\{i,j\},J\cup\{k,\ell\}}$};
        \node at (6.5,1.1) {$\scriptstyle\Delta_{I\cup\{j\},J\cup\{k\}}$};
        \node at (3.5,1.05) {$\scriptstyle\Delta_{I\cup\{j\},J\cup\{\ell\}}$};
        \node at (0.5,1) {$\scriptstyle\Delta_{I\cup\{i\},J\cup\{\ell\}}$};
        \node at (3.5,-0.5) {$\scriptstyle\Delta_{I,J}$};
        \node at (3.5,-1) {};
      \end{tikzpicture}
      \caption{Chamber minors associated to a simple transposition of opposite letters in a double reduced word.}
      \label{fig:wiring exchange relation}
    \end{figure}

    In a double wiring diagram, we label the strands in $u$'s wiring diagram on the right by $1$ through $n$ starting from the bottom and label the strands in $v$'s wiring diagram similarly on the left.  In this way, each chamber of the double wiring diagram determines a minor with row set given by the labels of strands in $v$'s wiring diagram lying below the chamber and column set given by the labels of strands in $u$'s wiring diagram lying below the chamber (see Figure~\ref{fig:double wiring diagram}).  Now we may describe a collection of total positivity criteria for $GL_n^{u,v}$.
    \setcounter{subtheorem}{3}
    \begin{theorem}\label{th:total positivity criteria}\cite{FZ99}
      Each double wiring diagram for $(u,v)$ determines a total positivity criterion for $GL_n^{u,v}$, i.e. an element of $GL_n^{u,v}$ lies in $GL_n^{\ge0}$ if and only if all $n+\ell(u)+\ell(v)$ chamber minors determined by the double wiring diagram are positive.
    \end{theorem}
    The totally positive matrices are exactly the elements of $GL_n^{w_0,w_0}\cap GL_n^{\ge0}$.  Thus Theorem~\ref{th:total positivity criteria} provides many total positivity criteria for $GL_n$.
    \setcounter{subexercise}{4}
    \begin{exercise}
      Use double wiring diagrams to find a total positivity criterion for $GL_3$.  Does the total positivity criterion you found in Exercise~\ref{exer:gl3 positivity} come from a double wiring diagram?
    \end{exercise}
    To finish the section and connect to the theory of cluster algebras to be presented in the next section we make the following observation.  Suppose a double reduced word for $(u,v)$ contains neighboring letters, one belonging to a reduced word for $u$ and one belonging to a reduced word for $v$.  In Figure~\ref{fig:wiring exchange relation}, we show the local effect on a double wiring diagram when these letters are interchanged. 

    Observe that under such an exchange we have the following exchange relation analogous to \eqref{eq:short plucker special}:
    \begin{equation}\label{eq:short Plucker}
      \Delta_{I\cup\{i\},J\cup\{k\}}\Delta_{I\cup\{j\},J\cup\{\ell\}}=\Delta_{I,J}\Delta_{I\cup\{i,j\},J\cup\{k,\ell\}}+\Delta_{I\cup\{i\},J\cup\{\ell\}}\Delta_{I\cup\{j\},J\cup\{k\}}.
    \end{equation}
    Moreover, observe that any two collections of chamber minors for double reduced words of $(u,v)$ can be obtained from each other by a sequence of these exchanges.
    \begin{exercise}
      Prove equation \eqref{eq:short Plucker}.
    \end{exercise}
		
		\begin{exercise}
		Another local transformation of a double wiring diagram is given by braid moves, either within $u$ or within $v$.  As an example, suppose a doubled word for $u,v$ begins $1,2,1$ (i.e.\ $u$ begins $1,2,1$ and these all occur before the first letter of $v$ in the doubled word).  Consider the double wiring diagrams both for this word and the one obtained by replacing the first three letters with $2,1,2$.  Show that the two corresponding collections of chamber minors differ in one element.  Find and prove an exchange relation that describes this transformation.
		\end{exercise}



\section{A Unifying Concept: Cluster Algebras}\label{sec:cluster_algebras}
	We now define cluster algebras of geometric type and $Y$-patterns, with a focus on the underlying dynamics of seed mutations. All definitions are due to S. Fomin and A. Zelevinsky and are drawn from \cite{FZ02} and \cite{FZ07}.
	
	\subsection{Basic definitions}
	Each cluster algebra is defined recursively from some initial data called a seed.  A seed consists of a cluster, which was informally described at the beginning of Section \ref{sec:motivation}, together with some combinatorial data that encode the mutations that can be performed.

	\begin{definition}
		Let $\cF$ be a purely transcendental field extension of $\CC$.  A \emph{seed}\index{seed} is a pair $(\x,B)$ where $\x = (x_1,\ldots, x_n)$ is an $n$-tuple of elements forming a transcendence basis of $\cF$ over $\CC$ and $B$ is a skew-symmetric integer $n \times n$ matrix.  The collection $\x$ is called the \emph{cluster}\index{cluster} and the matrix $B$ is called the \emph{exchange matrix}\index{exchange matrix}.
	\end{definition}
  
	The following employs the notation $[a]_+ = \max(a,0)$.
	\begin{definition}
		Given a seed ($\x,B$) and an integer $k=1,2,\ldots, n$ the \emph{seed mutation}\index{seed!1@mutation} $\mu_k$ in direction $k$ produces a new seed $\mu_k(\x,B) = (\x',B')$ where $\x' = (x_1,\ldots, x_{k-1}, x_k', x_{k+1},\ldots, x_n)$ with 
		\begin{equation} \label{eq:exchange relation}
		x_k' = \frac{\prod_{b_{ik}>0} x_i^{b_{ik}} + \prod_{b_{ik}<0} x_i^{-b_{ik}}}{x_k}
		\end{equation}
		and $B'$ is defined by
		\begin{equation}\label{eq:matrix mutation}
      b'_{ij}=\begin{cases}
                 -b_{ij} & \text{if $i=k$ or $j=k$;}\\
                 b_{ij}+[b_{ik}]_+[b_{kj}]_+-[-b_{ik}]_+[-b_{kj}]_+ & \text{otherwise.}
               \end{cases}
    \end{equation}
	\end{definition}
		
	In words, the mutation $\mu_k$ has the following effects:
	\begin{enumerate}
	\item $x_k$ changes to $x_k'$ satisfying $x_kx_k' = (\,\textrm{binomial in the other } x_i$);
	\item the entries $b_{ij}$ of $B$ away from row and column $k$ increase (resp. decrease) by $b_{ik}b_{kj}$ if $b_{ik}$ and $b_{kj}$ are both positive (resp. both negative);
	\item the entries of the $k^{th}$ row and the $k^{th}$ column of $B$ are negated.
	\end{enumerate}
	
	\begin{lemma} \label{lem:mutate}
	Let $(\x,B)$ be a seed in $\cF$.  For $k=1,\ldots, n$, the following hold:
	\begin{enumerate}
		\item $\mu_k(\x,B)$ is also a seed in $\cF$;
		\item the seed mutation $\mu_k$ is involutive, i.e. $\mu_k(\mu_k(\x,B)) = (\x,B)$.
	\end{enumerate}
	\end{lemma}
	
	\begin{definition}
		Fix an ambient field $\cF$ and an initial seed $(\bx, B)$.  The entries of the clusters of all seeds reachable from this one by a sequence of mutations are called the \emph{cluster variables}\index{cluster!1@variable}.  The \emph{cluster algebra}\index{cluster!2@algebra} associated with the initial seed is the subalgebra $\cA:=\cA(\bx,B)$ of $\cF$ generated by the set of all cluster variables.
	\end{definition}
	
  \begin{example}\label{example:type A2}
		Let $\x = (x_1,x_2)$ be the initial cluster with initial exchange matrix
		\begin{displaymath}
			B = \left[ \begin{array}{cc}
			0 & 1 \\
			-1 & 0 \\
			\end{array}
			\right].
		\end{displaymath}
		Then $\mu_1(\x,B) = ((x_1',x_2), -B)$ where
		\begin{displaymath}
			x_1' = \frac{x_2+1}{x_1}.
		\end{displaymath}
		It is convenient to denote the new cluster variable $x_1'=x_3$.  Next $\mu_2((x_3,x_2),-B) = ((x_3,x_4),B)$ where
		\begin{displaymath}
			x_4 = \frac{x_3+1}{x_2} = \frac{\frac{x_2+1}{x_1}+1}{x_2} = \frac{x_1 + x_2 + 1}{x_1x_2}
		\end{displaymath}
		and  $\mu_2((x_3,x_4),B) = ((x_5,x_4),-B)$ where
		\begin{displaymath}
			x_5 = \frac{x_4+1}{x_3} = \frac{\frac{x_1 + x_2 + 1}{x_1x_2}+1}{\frac{x_2+1}{x_1}} = \frac{x_1+1}{x_2}.
		\end{displaymath}
		The mutation pattern seems clear, but remarkably the next variable $x_6 = (1+x_5)/x_4$ equals the first variable $x_1$.  In fact, the only distinct cluster variables that can be obtained are $x_1$ through $x_5$, so
		\begin{displaymath}
			\cA(\x,B) = \CC\left[x_1,x_2,\frac{x_2+1}{x_1},\frac{x_1 + x_2 + 1}{x_1x_2},\frac{x_1+1}{x_2}\right] \subseteq \CC[x_1^{\pm1},x_2^{\pm1}].
		\end{displaymath}
		Note the generating set is not minimal, in fact, any four of the given elements generate $\cA(\x,B)$.
  \end{example}
	
	In the preceding example, there were only finitely many cluster variables.  As we will see, this is typically not the case.  However, a more subtle feature does hold in general, namely that each of the cluster variables is a Laurent polynomial in the variables of the initial seed.  
	\begin{theorem}\label{th:Laurent phenomenon}
		For any initial seed $(\x,B)$ with $\x = (x_1,\ldots, x_n)$, the associated cluster algebra lies in the Laurent polynomial ring
		\begin{displaymath}
			\cA(\x,B) \subseteq \CC[x_1, x_1^{-1}, \ldots, x_n, x_n^{-1}].
		\end{displaymath}
		In particular, every cluster variable can be expressed as a Laurent polynomial in $x_1,\ldots, x_n$.
	\end{theorem}
	Section \ref{sec:basic_theorems} will go into more detail on this result, called the Laurent phenomenon\index{Laurent phenomenon}, as well as several more main theorems on cluster algebras.
	
	\subsection{Increased generality}
	We generalize the previous definitions in two ways, first by allowing more general exchange matrices and then by allowing for certain coefficients in the exchange relations.
	
	\begin{definition}
		An $n \times n$ integer matrix $B$ is \emph{skew-symmetrizable}\index{skew-symmetrizable matrix} if there is a diagonal matrix $D$ with positive integer diagonal entries such that $DB$ is skew-symmetric.
	\end{definition}
	
	As an example, a $2 \times 2$ matrix $B$ is skew-symmetrizable if and only if
	\begin{displaymath}
		B = \left[ \begin{array}{cc} 	0 & b \\	-c & 0 \\	\end{array}	\right]
		\textrm{ or } B = \left[ \begin{array}{cc} 	0 & -b \\	c & 0 \\	\end{array}	\right]
	\end{displaymath}
	for positive integers $b$ and $c$.  In either case, a possible symmetrizing matrix is
	\begin{displaymath}
		D = \left[ \begin{array}{cc} 	c & 0 \\	0 & b \\	\end{array}	\right]
	\end{displaymath}
	The \emph{rank} of a cluster algebra is the number $n$ of elements in each cluster, so the exchange matrices $B$ just described give rise to all possible rank $2$ cluster algebras.
	
	The second generalization is to allow so-called \emph{frozen variables}\index{frozen variable} that never mutate, but which play a part in the exchange relations for the cluster variables. An \emph{extended cluster}, by convention, is typically written $\x = (x_1,\ldots, x_n, x_{n+1},\ldots x_m)$ where $x_1,\ldots, x_n$ are the cluster variable and $x_{n+1},\ldots, x_m$ are the frozen variables.  An \emph{extended exchange matrix} is an $m \times n$ integer matrix $\tilde{B}$ with the property that its upper $n \times n$ submatrix is skew-symmetrizable.
	
	\begin{definition}
	Fix a seed $(\x, \tilde B)$ with $\x=(x_1,\ldots, x_m)$ an extended cluster and $\tilde{B}$ an extended exchange matrix.  For an integer $k = 1,2,\ldots, n$, the \emph{seed mutation} $\mu_k$ in direction $k$ produces a new seed $\mu_k(\x,\tilde{B}) = (\x',\tilde{B}')$ with $\x' = (x_1,\ldots, x_{k-1}, x_k', x_{k+1},\ldots, x_m)$.  The formulas for $x_k'$ and the entries $\tilde{b}'_{ij}$ of $B'$ are the same as in \eqref{eq:exchange relation} and \eqref{eq:matrix mutation}, where the products in \eqref{eq:exchange relation} now range from $1$ to $m$ instead of from $1$ to $n$.
	\end{definition}
	
	\begin{lemma}
	Lemma \ref{lem:mutate} holds in this generalized setting.  Moreover, the skew-symmetrizing matrix $D$ is unchanged by mutation.
	\end{lemma}
	
	Given a seed $(\x, \tilde{B})$ as above, the corresponding cluster algebra\index{cluster!2@algebra} $\cA(\x, \tilde{B})$ is defined to be the subalgebra of $\CC(x_1,\ldots, x_m)$ generated by all cluster variables reachable from this seed together with the frozen variables (which appear in every seed).  Cluster algebras at this level of generality are referred to as cluster algebras of \emph{geometric type}.  In comparison, the earlier definition given was of \emph{skew-symmetric} cluster algebras with \emph{trivial coefficients}.
	
	\begin{example} \label{example:rank2}
		Consider a rank $2$ cluster algebras with trivial coefficients (so $n=m=2$).  The initial exchange matrix can be taken to be
		\begin{displaymath}
			B = \left[ \begin{array}{cc} 	0 & b \\	-c & 0 \\	\end{array}	\right].
		\end{displaymath}
		with $b,c$ positive integers.  Then $\mu_1(B) = \mu_2(B) = -B$ and $\mu_1(-B) = \mu_2(-B) = B$, so as in Example \ref{example:type A2} (the case where $b=c=1$) the cluster variables can be identified as $\{x_i : i \in \mathbb{Z}\}$ (possibly with redundant labels) and the mutations given by 
		\begin{displaymath}
		\ldots \stackrel{\mu_1}{\longleftrightarrow} ((x_1,x_0),-B) \stackrel{\mu_2}{\longleftrightarrow} ((x_1,x_2),B)) \stackrel{\mu_1}{\longleftrightarrow} ((x_3,x_2),-B)) \stackrel{\mu_2}{\longleftrightarrow} ((x_3,x_4),B)) \stackrel{\mu_1}{\longleftrightarrow} \ldots
		\end{displaymath}
		The exchange relations are
		\begin{displaymath}
		x_{k-1}x_{k+1}=\begin{cases}x_k^b+1 & \text{if $k$ is odd;}\\ x_k^c+1 & \text{if $k$ is even.}\end{cases}
		\end{displaymath}
		
		\begin{subexercise}\label{exercise:rank 2}\mbox{}
            \begin{enumerate}
                \item Compute all cluster variables and coefficient variables for cluster algebras associated to $B=\left[\begin{array}{cc} 0 & b\\ -c & 0\end{array}\right]$ with $b,c\in\ZZ_{>0}$ and $bc\le 3$.
                \item Justify why attempting such a calculation is futile for $bc\ge4$.  Hint: consider the degrees appearing in the denominators of the cluster variables.
    		\end{enumerate}
  		\end{subexercise}
	
		\begin{subexercise}\label{exercise:rank 2 Laurent phenomenon}
    		Prove the Laurent phenomenon for rank 2 cluster algebras.  Hint: start by showing that $x_4\in\ZZ[x_0,x_1,x_2,x_3]$ by considering the monomial $x_0x_3^b$.
  		\end{subexercise}
	\end{example}
	
	\begin{example}
		Let $n=3$ and $m=9$.  Denote the cluster 
		\begin{displaymath}
			\x = (\Delta_{13},\Delta_{14},\Delta_{15}, \Delta_{12},\Delta_{23},\Delta_{34},\Delta_{45},\Delta_{56},\Delta_{16})
		\end{displaymath}
		where the last $m-n=6$ variables are frozen.  Let
		\begin{displaymath}
			\tilde{B} = \left[\begin{array}{ccc}
			0 & 1 & 0 \\
			-1 & 0 & 1 \\
			0 & -1 & 0 \\
			1 & 0 & 0 \\
			-1 & 0 & 0 \\
			1 & -1 & 0 \\
			0 & 1 & -1 \\
			0 & 0 & 1 \\
			0 & 0 & -1 \\
			\end{array}
			\right].
		\end{displaymath}
		Let $\cA = \cA(\x,\tilde{B})$.  Interpreting the initial variables as the indicated Pl\"ucker coordinates on $Gr_{2,6}$ we obtain a cluster algebra in the ring of rational functions on $Gr_{2,6}$.
		\begin{subexercise}
            Prove the following statements about the cluster algebra $\cA$ just defined:
			\begin{itemize}
			\item The cluster variables of $\cA$ are precisely the Pl\"ucker coordinates $\Delta_{ij}$ with $1 \leq i < j \leq n$.
			\item The clusters of $\cA$ are precisely the collections $\Delta_T$ (see Section \ref{sub:Gr2n}) for $T$ a triangulation of a hexagon. 
			\end{itemize}
		\end{subexercise}
		For instance, the new cluster variable produced by applying the mutation $\mu_1$ to the initial seed is
		\begin{displaymath}
		\frac{\Delta_{12}\Delta_{34} + \Delta_{14}\Delta_{23}}{\Delta_{13}}
		\end{displaymath}
		which equals $\Delta_{24}$.
	\end{example}
	
	There is an alternate formulation of cluster algebras of geometric type which focuses on the roles played by the frozen variables $x_{n+1},\ldots, x_m$ in the exchange relations, rather than on the variables themselves.  Let $(\x, \tilde{B})$ be an extended seed.  For $k=1,\ldots, n$ let 
	\begin{displaymath}
	y_k = \prod_{i=n+1}^m x_i^{b_{ik}}.
	\end{displaymath}
	Define an operation $\oplus$ (called \emph{auxiliary addition}) on Laurent monomials by
	\begin{displaymath}
	\prod_{i=n+1}^m x_i^{e_i} \oplus \prod_{i=n+1}^m x_i^{f_i} = \prod_{i=n+1}^m x_i^{\min(e_i,f_i)}.
	\end{displaymath}
	Using this operation, we can extract positive and negative exponents as follows
	\begin{displaymath}
	\frac{y_k}{1 \oplus y_k} = \prod_{\substack{i=n+1,\ldots, m\\ b_{ik > 0}}} x_i^{b_{ik}}
	\quad \quad \quad \frac{1}{1 \oplus y_k} = \prod_{\substack{i=n+1,\ldots, m\\ b_{ik < 0}}} x_i^{-b_{ik}}.
	\end{displaymath}
	The change in going from the original exchange relation \eqref{eq:exchange relation} to the one in geometric type can then be summarized by saying that the two terms of the binomial are enriched with coefficients $y_k/(1 \oplus y_k)$ and $1/(1 \oplus y_k)$.  We now write
	\begin{displaymath} 
		x_k' = \frac{y_k\prod_{b_{ik}>0} x_i^{b_{ik}} + \prod_{b_{ik}<0} x_i^{-b_{ik}}}{(1 \oplus y_k)x_k}
	\end{displaymath}
	where the products range over $i$ from $1$ to $n$.
	
	\subsection{$Y$-patterns}
	There is an alternate version of seeds and mutations, closely related to the previous one, which itself arises in many applications.  The resulting structures are called \emph{$Y$-patterns}\index{Y-pattern}.
	
	\begin{definition}
		A \emph{$Y$-seed} is a pair $(\y,B)$ consisting of an $n$-tuple $\y=(y_1,\ldots, y_n)$ and an $n\times n$ skew-symmetric matrix $B$.  For $k=1,\ldots, n$, the \emph{$Y$-seed mutation} $\mu_k$ is defined by
		\begin{displaymath}
			\mu_k((y_1,\ldots, y_n),B) = ((y_1',\ldots, y_n'),B')
		\end{displaymath}
		for $B'$ as defined in \eqref{eq:matrix mutation} and 
		\begin{equation}\label{eq:y mutation}
		y_j' = \begin{cases}
		y_k^{-1} & \text{if $j = k$;} \\
		y_jy_k^{[-b_{jk}]_+}(1+y_k)^{b_{jk}} & \text{if $j \neq k$.} \\
		\end{cases}
		\end{equation}
	\end{definition}
	
	In words, the $Y$-seed mutation $\mu_k$ has the following effects:
	\begin{enumerate}
	\item For each $j \neq k$, $y_j$ is multiplied by $(1+y_k)^{b_{jk}}$ if $b_{jk}>0$ and by $(1+y_k^{-1})^{b_{jk}}$ if $b_{jk}<0$,
	\item $y_k$ is inverted,
	\item $B$ is changed in the same way as for ordinary seed mutations.
	\end{enumerate}
	The $Y$-dynamics also go by the name \emph{coefficient dynamics}\index{coefficient dynamics} because the coefficients in a geometric type cluster algebra evolve in this manner with respect to auxiliary addition $\oplus$ and ordinary multiplication of Laurent monomials.
	
	Another connection between cluster algebra and $Y$-pattern dynamics comes by way of a certain Laurent monomial change of variables.  Let $\tilde{B}$ be an $m \times n$ extended exchange matrix with $B$ its upper $n \times n$ submatrix.  Given a seed $((x_1,\ldots, x_m),\tilde{B})$, define an associated $Y$-seed $((\hat{y}_1,\ldots, \hat{y}_n),B)$ by
	\begin{displaymath}
	\hat{y}_j = \prod_{i=1}^m x_i^{b_{ij}}
	\end{displaymath}
	
	\begin{proposition}\label{prop:yhat mutations}
	Fix $k=1,\ldots, n$ and suppose $\mu_k(\x,\tilde{B}) = (\x',\tilde{B}')$.  Define $(\hat{y}_1,\ldots, \hat{y}_n)$ from $(\x,\tilde{B})$ and $(\hat{y}_1',\ldots, \hat{y}_n')$ from $(\x',\tilde{B}')$ as above.  Then 
	\begin{displaymath}
	\mu_k((\hat{y}_1,\ldots, \hat{y}_n), B) = ((\hat{y}_1',\ldots, \hat{y}_n'), B').
	\end{displaymath}
	\end{proposition}

  \begin{exercise}
    Prove Proposition~\ref{prop:yhat mutations}.
  \end{exercise}

\section{Foundational Results}\label{sec:basic_theorems}
  As mentioned in the introduction cluster algebras have found applications in a surprising array of mathematical disciplines.  Much of their ubiquity comes from a number of remarkable theorems which we now explain.  

  The first of these results is the aforementioned Laurent phenomenon (Theorem \ref{th:Laurent phenomenon}).  The Laurent property originates from the study of several classical recurrences that seem as though they should produce rational numbers, but turn out to give only integers under certain initial conditions \cite{FZ02b}.  As with several of the results in this section, we give a complete proof in the case of rank 2 cluster algebras.
	
	Recalling Example~\ref{example:rank2}, consider the rank 2 cluster algebra with initial exchange matrix $B=\left[\begin{array}{cc}0 & b\\ -c & 0\end{array}\right]$  and cluster variables $x_k$, $k\in\ZZ$ which satisfy the recursion
  \[x_{k-1}x_{k+1}=\begin{cases}x_k^b+1 & \text{if $k$ is odd;}\\ x_k^c+1 & \text{if $k$ is even.}\end{cases}\]
  \begin{theorem}
    For any $m\in\ZZ$, each $x_k$ is an element of $\ZZ[x_{m-1},x_m,x_{m+1},x_{m+2}]$.
  \end{theorem}
  \begin{remark}
    The ring $\ZZ[x_{m-1},x_m,x_{m+1},x_{m+2}]$ is called the \emph{lower bound} of $\cA(\bx,B)$ at the cluster $\{x_m,x_{m+1}\}$.  We will give more details below and see that this result is a special case of Theorem~\ref{th:lower bounds}.
  \end{remark}
  \begin{proof}
		Assume without loss of generality that $m$ is odd (otherwise interchange $b$ and $c$ in the calculation below).  Then we may compute
    \begin{align*}
      x_{m-1}x_{m+2}^b
      &=\frac{x_m^b+1}{x_{m+1}}x_{m+2}^b\\
      &=\frac{x_m^bx_{m+2}^b-1}{x_{m+1}}+\frac{x_{m+2}^b+1}{x_{m+1}}\\
      &=\frac{(x_{m+1}^c+1)^b-1}{x_{m+1}}+x_{m+3}.
    \end{align*}
    In particular, we see that $x_{m+3}\in\ZZ[x_{m-1},x_m,x_{m+1},x_{m+2}]$.  A similar calculation shows $x_{m-2}\in\ZZ[x_{m-1},x_m,x_{m+1},x_{m+2}]$ and a simple induction argument of ``shifting the viewing window'' establishes the result.
  \end{proof}
	Having shown $x_k$ to be a polynomial in $x_{m-1},x_m,x_{m+1},x_{m+2}$, each of which is directly seen to be a Laurent polynomial in $x_m, x_{m+1}$, the rank 2 Laurent phenomenon follows.  As such, the above is a solution to Exercise \ref{exercise:rank 2 Laurent phenomenon}. With more work one may establish the following result which can be leveraged to prove the Laurent phenomenon in general.
  \begin{theorem}
    Assume $b,c\ne0$.  For any $m\in\ZZ$, the cluster algebra $\cA(\bx,B)$ is equal to $\bigcap_{k\in\ZZ}\ZZ[x_k^{\pm1},x_{k+1}^{\pm1}]=\bigcap_{k=m-1}^{m+1}\ZZ[x_k^{\pm1},x_{k+1}^{\pm1}]$.
  \end{theorem}
  The proof of this result is somewhat technical and seems unlikely to be very informative so we will omit the details.
  \begin{remark}
    If one allows extended exchange matrices and imposes an additional coprimality assumption (see Definition~\ref{def:coprime}), then the condition $b,c\ne0$ may be dropped.  The ring $\bigcap_{k=m-1}^{m+1}\ZZ[x_k^{\pm1},x_{k+1}^{\pm1}]$ is called the \emph{upper bound} of $\cA(\bx,B)$ at the cluster $\{x_m.x_{m+1}\}$ and the ring $\bigcap_{k\in\ZZ}\ZZ[x_k^{\pm1},x_{k+1}^{\pm1}]$ is called the \emph{upper cluster algebra}.  We will give more details below, at which point this result will be a special case of Theorems~\ref{th:lower bounds} and~\ref{th:acyclic upper bounds}.
  \end{remark}

	The next result is the classification of cluster algebras of \emph{finite type}, i.e.\ those that have only finitely many cluster variables.  As we have seen in Sections~\ref{sub:Gr2n} and~\ref{sub:total positivity}, cluster algebras have roots in classical Lie theory.  A large number of objects in this realm are classified by Dynkin diagrams and finite type cluster algebras are no exception.  Given an extended exchange matrix $\tilde B$ with principal square submatrix $B=(b_{ij})_{i,j=1}^n$, write $A:=A_{\tilde B}=(a_{ij})_{i,j=1}^n$ for the \emph{Cartan companion}\index{Cartan companion} of $\tilde B$ given by 
  \[a_{ij}=\begin{cases}2 & \text{if $i=j$;}\\ -|b_{ij}| & \text{otherwise.}\end{cases}\]
  Note that the Cartan companion only depends on the principal part of $\tilde B$.
  \begin{theorem}\cite{FZ03}
    A cluster algebra $\cA(\x,\tilde B)$ is of finite type if and only if there exists an extended exchange matrix mutation equivalent to $\tilde B$ whose Cartan companion is a finite-type Cartan matrix.
  \end{theorem}
  \begin{example}
		The irreducible finite type $2 \times 2$ Cartan matrices are (up to transposition) given by
		\begin{align*}
		&A_2 : \left[\begin{array}{cc} 2 & -1\\ -1 & 2\end{array}\right] \\
		&B_2 : \left[\begin{array}{cc} 2 & -2\\ -1 & 2\end{array}\right] \\
		&G_2 : \left[\begin{array}{cc} 2 & -3\\ -1 & 2\end{array}\right]
		\end{align*}
		or put another way, they are the matrices $\left[\begin{array}{cc} 2 & -b\\ -c & 2\end{array}\right]$ with $bc \leq 3$.  This is precisely the Cartan counterpart $A_B$ of the exchange matrix $B=\left[\begin{array}{cc} 0 & b\\ -c & 0\end{array}\right]$.  Therefore, the rank 2 cluster algebra with this exchange matrix, or any cluster algebra obtained from it by adding coefficients, has finite type if and only if $bc \leq 3$.  This was the content of Exercise \ref{exercise:rank 2}.
		\end{example}
		
		We now provide an outline to the solution of Exercise \ref{exercise:rank 2}.  Let $B=\left[\begin{array}{cc} 0 & b\\ -c & 0\end{array}\right]$ and let $x_k, k \in \mathbb{Z}$ be the cluster variables of this cluster algebra.  We need some more notation.  For $k \neq 1,2$, let $d^{(k)}=(d_1^{(k)},d_2^{(k)})$ denote the exponents of the monomial $x_1^{d^{(k)}_1}x_2^{d^{(k)}_2}$ appearing in the denominator of $x_k$ when expressed as a Laurent polynomial in $x_1, x_2$.  Let $U_\ell(t)$ denote the (normalized) Chebyshev polynomial\index{Chebyshev polynomial} (of the second kind) defined recursively by 
    \[U_{\ell+1}(t)=tU_\ell(t)-U_{\ell-1}(t)\qquad U_1(t)=1\qquad U_0(t)=0\]
    and for $j\in\{1,2\}$ define 
    \[u_{\ell,j}=\begin{cases}U_\ell(\sqrt{bc}) & \text{if $\ell$ is odd;}\\\sqrt{b/c}\,U_\ell(\sqrt{bc}) & \text{if $\ell$ is even and $j=1$;}\\\sqrt{c/b}\,U_\ell(\sqrt{bc}) & \text{if $\ell$ is even and $j=2$.}\\\end{cases}\]
    Note that $U_\ell(t)$ is an odd (resp. even) function when $\ell$ is even (resp. odd), so each $u_{\ell,j}$ is a polynomial in $b$ and $c$.  Moreover, one may easily show that each $u_{\ell,j}$, $\ell\ge1$ is a positive integer when $bc\ge4$.  
	
		The rank 2 finite type classification now boils down to verifying the following:
		\begin{itemize}
		\item The denominator vector of $x_k$ for $k\in\ZZ\setminus\{1,2\}$ is
    \begin{equation}\label{eq:rank 2 denominators}
      d^{(k)}=\begin{cases}(u_{k-2,1},u_{k-3,2}) & \text{if $k\ge3$;}\\(u_{-k,1},u_{-k+1,2}) & \text{if $k\le0$.}\end{cases}
    \end{equation}
    \item $U_n\big(2\cos(\theta)\big)=\frac{\sin(n\theta)}{\sin(\theta)}$ from which it follows that there are infinitely many $d^{(k)}$ if and only if $bc\ge4$.  In particular, there are infinitely many cluster variables if $bc \geq 4$.
		\item By explicit calculation, there are only finitely many cluster variables in each case where $bc \leq 3$.
    \end{itemize}
 
    
  More generally, if $\tilde{B}$ is any $m\times n$ extended exchange matrix and $|b_{ij}b_{ji}| \geq 4$ for some $i,j=1,\ldots, n$, then iteratively mutating at $i$ and $j$ alone will produce infinitely many cluster variables.  Call an extended exchange matrix $\tilde B=(b_{ij})$ \emph{2-finite} if every matrix $\tilde B'=(b'_{ij})$ mutation equivalent to $\tilde B$ satisfies the condition $|b'_{ij}b'_{ji}|\le3$ for all $i,j$.  By the above reasoning, 2-finiteness is a necessary condition for a cluster algebra to be of finite type.  It is quite remarkable that it is sufficient as well.
  \begin{theorem}\cite{FZ03}
    A cluster algebra $\cA(\x,\tilde B)$ is of finite type if and only if $\tilde B$ is 2-finite.
  \end{theorem}

  The next result gives the positivity of the initial cluster Laurent expansions of all cluster variables.  This result was conjectured with the introduction of cluster algebras \cite{FZ02} and remained open for more than ten years.  
  \begin{theorem}\cite{LS15,GHKK14}
    Let $\cA(\bx,\tilde B)$ be a cluster algebra.  Any cluster variable of $\cA(\bx,\tilde B)$ is an element of $\ZZ_{\ge0}[x_1^{\pm1},\ldots,x_m^{\pm1}]$.
  \end{theorem}
  The proof in skew-symmetric types builds on a concrete combinatorial construction of cluster variables in rank 2 \cite{LS13,LLZ14}.  For the more general skew-symmetrizable cluster algebras the proof was given in \cite{GHKK14} using the theory of scattering diagrams built on connections to mirror symmetry.  We describe the combinatorial approach to the skew-symmetric case here, but first we need more notation.

  For $k\in\ZZ\setminus\{1,2\}$, denote by $R_k$ the rectangle in $\ZZ^2$ with corner vertices $(0,0)$ and $d^{(k)}$ from Equation~\ref{eq:rank 2 denominators}.  Write $D_k$ for the maximal Dyck path\index{maximal Dyck path} in $R_k$ beginning at $(0,0)$, taking East and North steps to end at the upper right corner of $R_k$ while never passing above the main diagonal, and such that the area below $D_k$ inside $R_k$ is maximized.
  \begin{example}
    For $b=c=3$, the maximal Dyck path $D_5$ is
    
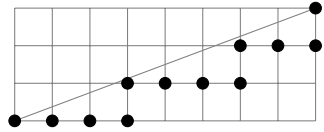
\begin{center}
    \begin{tikzpicture}
      \draw[step=0.5cm,color=gray] (0,0) grid (4.0,1.5);
      \draw[color=gray] (0,0) -- (4.0,1.5);
      \draw[fill=black] (0.0,0.0) circle (2.2pt);
      \draw[fill=black] (0.5,0.0) circle (2.2pt);
      \draw[fill=black] (1.0,0.0) circle (2.2pt);
      \draw[fill=black] (1.5,0.0) circle (2.2pt);
      \draw[fill=black] (1.5,0.5) circle (2.2pt);
      \draw[fill=black] (2.0,0.5) circle (2.2pt);
      \draw[fill=black] (2.5,0.5) circle (2.2pt);
      \draw[fill=black] (3.0,0.5) circle (2.2pt);
      \draw[fill=black] (3.0,1.0) circle (2.2pt);
      \draw[fill=black] (3.5,1.0) circle (2.2pt);
      \draw[fill=black] (4.0,1.0) circle (2.2pt);
      \draw[fill=black] (4.0,1.5) circle (2.2pt);
    \end{tikzpicture}
    \end{center}
  \end{example}
  For edges $e,e'\in D_k$, write $e<e'$ if $e$ precedes $e'$ along $D_k$, in other words if $e$ is closer to $(0,0)$.  In this case write $ee'$ for the subpath of $D_k$ beginning with the edge $e$ and ending with the edge $e'$.  Write $(ee')_H$ and $(ee')_V$ for the sets of horizontal and vertical edges in the path $ee'$ respectively.  Let $H,V\subset D_k$ denote the sets of horizontal and vertical edges of $D_k$.

  Call subsets $S_H\subset H$ and $S_V\subset V$ compatible\index{compatible pair} if for each $h\in S_H$ and $v\in S_V$ with $h<v$, there exists $e\in hv$ so that one of the following holds:
  \begin{itemize}
    \item $e\ne v$ and $|(he)_V|=c|(he)_H\cap S_H|$;
    \item $e\ne h$ and $|(ev)_H|=b|(ev)_V\cap S_V|$.
  \end{itemize}
  \begin{example}
    For $b=c=3$, the picture below on the left shows a compatible collection of edges of $D_5$ while the collection on the right is not compatible.
    \begin{center}
    \begin{tikzpicture}
      \draw[step=0.5cm,color=gray] (0,0) grid (4.0,1.5);
      \draw[color=gray] (0,0) -- (4.0,1.5);
      \draw[fill=black] (0.0,0.0) circle (2.2pt);
      \draw[fill=black] (0.5,0.0) circle (2.2pt);
      \draw[line width=2.5pt] (0.5,0.0) -- (1.0,0.0);
      \draw[fill=black] (1.0,0.0) circle (2.2pt);
      \draw[fill=black] (1.5,0.0) circle (2.2pt);
      \draw[fill=black] (1.5,0.5) circle (2.2pt);
      \draw[fill=black] (2.0,0.5) circle (2.2pt);
      \draw[fill=black] (2.5,0.5) circle (2.2pt);
      \draw[fill=black] (3.0,0.5) circle (2.2pt);
      \draw[line width=2.5pt] (3.0,0.5) -- (3.0,1.0);
      \draw[fill=black] (3.0,1.0) circle (2.2pt);
      \draw[fill=black] (3.5,1.0) circle (2.2pt);
      \draw[fill=black] (4.0,1.0) circle (2.2pt);
      \draw[line width=2.5pt] (4.0,1.0) -- (4.0,1.5);
      \draw[fill=black] (4.0,1.5) circle (2.2pt);
    \end{tikzpicture}
    \qquad
    \begin{tikzpicture}
      \draw[step=0.5cm,color=gray] (0,0) grid (4.0,1.5);
      \draw[color=gray] (0,0) -- (4.0,1.5);
      \draw[fill=black] (0.0,0.0) circle (2.2pt);
      \draw[fill=black] (0.5,0.0) circle (2.2pt);
      \draw[fill=black] (1.0,0.0) circle (2.2pt);
      \draw[line width=2.5pt] (1.0,0.0) -- (1.5,0.0);
      \draw[fill=black] (1.5,0.0) circle (2.2pt);
      \draw[fill=black] (1.5,0.5) circle (2.2pt);
      \draw[fill=black] (2.0,0.5) circle (2.2pt);
      \draw[fill=black] (2.5,0.5) circle (2.2pt);
      \draw[fill=black] (3.0,0.5) circle (2.2pt);
      \draw[line width=2.5pt] (3.0,0.5) -- (3.0,1.0);
      \draw[fill=black] (3.0,1.0) circle (2.2pt);
      \draw[fill=black] (3.5,1.0) circle (2.2pt);
      \draw[fill=black] (4.0,1.0) circle (2.2pt);
      \draw[line width=2.5pt] (4.0,1.0) -- (4.0,1.5);
      \draw[fill=black] (4.0,1.5) circle (2.2pt);
    \end{tikzpicture}
    \end{center}
    \begin{subexercise}
      Find all other compatible collections for $D_5$ in the case $b=c=3$.
    \end{subexercise}
  \end{example}
  \begin{theorem}\label{th:rank 2 positivity}\cite{LLZ14}
    For $k\in\ZZ\setminus\{1,2\}$, the cluster variable $x_k$ is given by
    \[x_k=\sum\limits_{(S_H,S_V)}x_1^{-d^{(k)}_1+b|S_V|}x_2^{-d^{(k)}_2+c|S_H|},\]
    where the sum ranges over all compatible collections of edges in the maximal Dyck path $D_k$.  In particular, each cluster variable is contained in $\ZZ_{\ge0}[x_1^{\pm1},x_2^{\pm1}]$.
  \end{theorem}

  \begin{remark}
    The idea of obtaining positivity for arbitrary cluster algebras goes as follows (see \cite{LS15}): any sequence of mutations can be viewed as a collection of rank 2 mutation sequences, so iteratively applying Theorem~\ref{th:rank 2 positivity} will (after a considerable amount of work) lead to a proof of positivity in general.
  \end{remark}

  One motivation for the discovery of the cluster algebra formalism was the desire to find a combinatorial construction of dual canonical basis elements for (quantum) algebraic groups.  The dual canonical basis of a semisimple algebraic group induces bases on the coordinate rings of any of its subvarieties.  As we saw in Example~\ref{sub:total positivity} there appears to be some kind of cluster algebra structure on the double Bruhat cells of a semisimple algebraic group, which one can hope will shed some light on the dual canonical basis.  In most cases the double Bruhat cell does not actually admit a cluster algebra structure but does admit a closely related structure.

  \begin{definition}
    Let $(\x,\tilde B)$ denote an extended seed.  Define the \emph{upper cluster algebra}\index{upper cluster algebra} $\overline{\cA}(\x,\tilde B)$ as the intersection of all Laurent rings associated to extended seeds mutation equivalent to $(\x,\tilde B)$.
  \end{definition}

  \begin{remark}
    The upper cluster algebra is exactly the collection of all rational functions in $\{x_1,\ldots,x_m\}$ which are expressed as Laurent polynomials in terms of any cluster.  Thus the Laurent phenomenon Theorem~\ref{th:Laurent phenomenon} establishes the inclusion $\cA(\x,\tilde B)\subset\overline{\cA}(\x,\tilde B)$, hence the name ``upper'' cluster algebra.
  \end{remark}

  \begin{theorem}\label{th:dbc is uca}\cite{BFZ05}
    The coordinate ring of any double Bruhat cell $GL_n^{u,v}$ is an upper cluster algebra such that each collection of chamber minors determines a cluster.
  \end{theorem}

  In some special cases, the cluster algebra coincides with its upper cluster algebra.  First we describe a necessary intermediate concept.
  \begin{definition}\label{def:coprime}
    Let $(\bx,\tilde B)$ be an extended seed and denote by $P_k:=x_kx_k'$ the binomial on the right hand side of the exchange relation \eqref{eq:exchange relation}.  We say that $(\bx,\tilde B)$ is \emph{coprime}\index{seed!2@coprime --} if the binomials $P_k$ are pairwise relatively prime.
  \end{definition}
  The \emph{upper bound}\index{upper bound algebra} $\cU(\bx',\tilde B)$ for $\cA(\bx,\tilde B)$ at a seed $(\bx',\tilde B')$ mutation equivalent to $(\bx,\tilde B)$ is the intersection of the Laurent ring generated by variables in $\bx'$ with all Laurent rings associated to seeds which can be obtained from $(\bx',\tilde B')$ by a single mutation.
  \begin{theorem}\cite{BFZ05}
    If $(\bx,\tilde B)$ and $(\bx',\tilde B')$ are mutation equivalent coprime seeds, then their upper bounds $\cU(\bx,\tilde B)$ and $\cU(\bx',\tilde B)$ coincide.  In particular, if every seed mutation equivalent to an extended seed $(\bx,\tilde B)$ is coprime, then $\overline{\cA}(\bx,\tilde B)=\cU(\bx',\tilde B')$ for any seed mutation equivalent to $(\bx,\tilde B)$.
  \end{theorem}

  \begin{definition}
    An exchange matrix $B=(b_{ij})$ is called \emph{acyclic}\index{exchange matrix!acyclic --} if there exists a permutation $\sigma$ so that $b_{\sigma_i\sigma_j}\le 0$ for $i<j$.
  \end{definition}
  This terminology can be easily understood in the case of a skew-symmetric exchange matrix; the matrix is acyclic exactly when its associated quiver (see Section \ref{sec:integrable_systems}) has no oriented cycles.
  \begin{theorem}\cite{BFZ05}\label{th:acyclic upper bounds}
    Let $\cA(\x,\tilde B)$ be a cluster algebra where $\tilde B$ is coprime and has acyclic principal part.  Then the cluster algebra $\cA(\x,\tilde B)$ coincides with its upper bound $\cU(\x,\tilde B)$.
  \end{theorem}

  Acyclicity also guarantees the existence of easily identifiable bases of a cluster algebra.  These are best understood by identifying the cluster algebra with another related algebra, its lower bound.
  \begin{definition}
    Let $(\x,\tilde B)$ be an extended seed and write $x'_k$ for the variable obtained by mutation in direction $k$.  The \emph{lower bound} of $\cA(\x,\tilde B)$ at $(\x,\tilde B)$ is the subalgebra $\cL(\x,\tilde B)=\CC[x_1,x'_1,\ldots,x_n,x'_n,x_{n+1},\ldots,x_m]$.
  \end{definition}

  \begin{theorem}\cite{BFZ05}\label{th:lower bounds}
    Let $\cA(\x,\tilde B)$ be a cluster algebra where $\tilde B$ is coprime and has acyclic principal part.  Then the cluster algebra $\cA(\x,\tilde B)$ coincides with its lower bound $\cL(\x,\tilde B)$.  Moreover, the collection of all \emph{standard monomials} in the variables $x_1,x'_1,\ldots,x_n,x'_n,x_{n+1},\ldots,x_m$, i.e. those which do not contain both $x_k$ and $x_k'$ for any $k$, forms a basis of the cluster algebra $\cA(\bx,\tilde B)$.
  \end{theorem}

  \begin{example}
    Acyclic cluster algebras are particularly nice but should not be taken as the general case: there exist cluster algebras where anything that can go wrong seems to go wrong.  The primary example of this is the \emph{Markov cluster algebra}\index{Markov!2@cluster algebra} with exchange matrix 
    \[B=\left[\begin{array}{ccc} 0 & 2 & -2\\ -2 & 0 & 2\\ 2 & -2 & 0\end{array}\right].\]
    
    Comparing with Example~\ref{sub:Markov triples} we see that the Markov cluster algebra has infinitely many cluster variables.  In this case, the following undesirable properties all hold:
    \begin{itemize}
      \item the standard monomials are linearly dependent;
      \item the Markov cluster algebra is not finitely generated and not Noetherian;
      \item the Markov cluster algebra (over $\CC$) contains non-prime cluster variables;
      \item the Markov cluster algebra does not equal its upper cluster algebra.
    \end{itemize}
  \end{example}

\section{Compatible Poisson Structures and Quantization}\label{sec:poisson_and_quantum}

  Many natural examples of cluster algebras carry the following additional structure. 
  \begin{definition}
    A \emph{Poisson algebra}\index{Poisson!algebra} $(A,\{\cdot,\cdot\})$ is an associative algebra $A$ equipped with an additional skew-symmetric bilinear operation, called the Poisson bracket\index{Poisson!bracket} (written as $\{x,y\}$ for $x,y\in A$), for which the following holds:
    \begin{itemize}
      \item given any $x\in A$, the endomorphism $\{x,\cdot\}:A\to A$, $y\mapsto\{x,y\}$ is a derivation with respect to both binary operations on $A$, i.e. for $x,y,z\in A$ we have
      \begin{align*}
        \tag{Leibnitz rule} \{x,yz\}&=\{x,y\}z+y\{x,z\}\\
        \tag{Jacobi identity}\{x,\{y,z\}\}&=\{\{x,y\},z\}+\{y,\{x,z\}\}.
      \end{align*}
    \end{itemize}
  \end{definition}

  In the case of cluster algebras, the Poisson brackets take an especially simple form when applied to cluster variables.  
  \begin{definition}
    Let $\cA$ be a cluster algebra.  A Poisson bracket $\{\cdot,\cdot\}$ on $\cA$ is compatible with the cluster algebra structure if every cluster $\bx$ of $\cA$ is \emph{log-canonical}\index{log-canonical basis} with respect to $\{\cdot,\cdot\}$, that is there exists a skew-symmetric matrix $\Omega_\bx=(\Omega_{ij})$ so that $\{x_i,x_j\}=\Omega_{ij}x_ix_j$ for all $i,j$.
  \end{definition}

  The existence of a compatible Poisson structure is quite restrictive and requires a particular relationship between the matrix $\Omega_\bx$ and the exchange matrix of a seed with cluster $\bx$. 
  \begin{lemma}
  \label{le:compatibility}
    Let $\cA(\bx,\tilde B)$ be a cluster algebra, where $\tilde B$ is $(m+n)\times n$, with compatible Poisson structure $\{\cdot,\cdot\}$.  Then writing $\Omega_\bx$ for the $(m+n) \times (m+n)$ matrix of Poisson bracket coefficients, we have $\tilde B^T\Omega_\bx=[D\, \boldsymbol{0}]$, where $D$ is a diagonal matrix which skew-symmetrizes $B$.
  \end{lemma}
  \begin{subexercise}
    Prove Lemma~\ref{le:compatibility} using the condition that each neighboring cluster to $\bx$ should also be log-canonical.
  \end{subexercise}
  
  In fact, these conditions can only hold in the following situation, in which case all possible compatible Poisson structures can be easily understood. 
  \begin{theorem}\cite{GSV03} \label{thm:Poisson}
		If an extended exchange matrix $\tilde B$ has full rank then a compatible Poisson structure exists for the associated cluster algebra.  Moreover, the collection of all such Poisson structures is parametrized by an affine space of dimension $\rho(B)+{m\choose 2}$, where $\rho(B)$ denotes the number of connected components of the quiver associated to the principal submatrix $B$ of $\tilde B$.
  \end{theorem}
  The primary motivating example for compatible Poisson structures comes in the form of the double Bruhat cells of Section~\ref{sub:total positivity}.  To introduce the relevant Poisson structure and place the double Bruhat cells in the proper context we need to introduce more notation.

  A smooth manifold $M$ together with a Poisson structure on its algebra $\cO(M)$ of regular functions is called a \emph{Poisson manifold}\index{Poisson!manifold}.  The best understood examples of such structures comes from the theory of symplectic geometry.
  \begin{example}
    A symplectic manifold\index{symplectic!1@manifold} $(M,\omega)$ is a smooth even-dimensional manifold $M^{2n}$ together with a 2-form $\omega\in H^2(M)$ that is non-degenerate, i.e. $\omega^n\ne0$ is a volume form.  The symplectic structure $\omega$ provides a natural association of a vector field $\xi_f$ to a function $f\in\cO(M)$ via the contraction formula $\iota_{\xi_f}\omega=-df$.  The algebra of smooth functions on the symplectic manifold $M$ is then naturally a Poisson algebra via $\{f,g\}=\omega(\xi_f,\xi_g)$ for $f,g\in\cO(M)$.
  \end{example}  
  A smooth map $\varphi:M\to N$ between Poisson manifolds is a \emph{Poisson morphism} if it induces a Poisson morphism $\varphi^*:\cO(N)\to\cO(M)$ on their algebras of smooth functions, i.e. if $\{f \circ \varphi,g \circ \varphi\}_M=\{f,g\}_N\circ\varphi$ for all $f,g\in \cO(N)$.  Given two Poisson manifolds $M$ and $N$, there is a natural Poisson structure on $M\times N$ given by
  \[\{f,g\}(x,y)=\{f(\cdot,y),g(\cdot,y)\}_M(x)+\{f(x,\cdot),g(x,\cdot)\}_N(y)\]
  for $f,g\in\cO(M\times N)$.
  \begin{example}
    The general linear group $GL_n$ is a Poisson manifold where the Poisson bracket on matrix entries is given by
    \begin{equation}
      \{x_{ij},x_{k\ell}\}=\frac{1}{2}\big(\sgn(k-i)+\sgn(\ell-j)\big)x_{i\ell}x_{kj}.
    \end{equation}
    In fact, $GL_n$ has the structure of a Poisson-Lie group, that is the multiplication map $GL_n\times GL_n\to GL_n$ is a Poisson morphism when $GL_n\times GL_n$ is given the product Poisson structure.
  \end{example}
  With this we may describe the Poisson structure on the double Bruhat cells.
  \begin{theorem}\cite{GSV10}
    The cluster algebra structure on $GL_n^{u,v}$ from Theorem~\ref{th:dbc is uca} is compatible with the restriction of the Poisson structure from $GL_n$.
  \end{theorem}

  The double Bruhat cells actually have a much closer connection to the Poisson geometry of $GL_n$.  To describe it we need to introduce a few additional concepts.

  The existence of a Poisson structure on a manifold allows one to construct vector fields associated to functions in a similar manner to the symplectic case.
  \begin{lemma}\label{le:Hamiltonian vector fields}
    Let $(M,\{\cdot,\cdot\})$ be a Poisson manifold.  For any $h\in\cO(M)$, there exists a vector field $\xi_h$ on $M$ such that for any $f\in\cO(M)$ we have $\xi_h(f)=\{f,h\}$.
  \end{lemma}
  The vector field $\xi_h$ from Lemma~\ref{le:Hamiltonian vector fields} is simply the vector field canonically associated to the derivation $\{\cdot,h\}$ and is called a \emph{Hamiltonian vector field}\index{Hamiltonian vector field}.  Flowing along Hamiltonian vector fields defines an equivalence relation on a Poisson manifold as follows: 
  \begin{quotation}
    $p\sim q$ if there exists a piecewise smooth curve connecting $p$ and $q$ where each smooth segment is the trajectory of a Hamiltonian vector field.
  \end{quotation}
  The equivalence class $M_p$ containing $p\in M$ is a connected submanifold of $M$ and by definition its tangent space $T_qM_p\subset T_qM$ is spanned by Hamiltonian tangent vectors.  We may view $M_p$ as a symplectic manifold whose 2-form is given by $\omega(\xi_f,\xi_g)=\{f,g\}$ and thus $M_p$ is called the \emph{symplectic leaf}\index{symplectic!2@leaf} of $M$ through $p$.
  \begin{example}
    The symplectic leaves of $GL_n$ have been explicitly computed (e.g. see \cite{KZ02} and references therein).  In particular, it is known that each double Bruhat cell in $GL_n$ is foliated by symplectic leaves of $GL_n$.  As an example (c.f. Example~\ref{example:big cell}) the following set is a symplectic leaf living in $GL_n^{w_0,w_0}$:
    \[\{M\in GL_n:\Delta_{[1,i],[n+1-i,n]}(M)=\Delta_{[i+1,n],[1,n-i]}(M)\ne0\text{ for all $i$}\}.\]
  \end{example}
  \begin{exercise}
    Find a matrix $M\in GL_2$ so that the symplectic leaf above is the symplectic leaf of $GL_2$ through $M$.  Do the same for $GL_3$.
  \end{exercise}

  A final reason to care about Poisson structures on a cluster algebra is that they provide a canonical quantization of the cluster algebra structure.  Rather than giving precise definitions, we indicate here the intuition motivating the definition and refer the reader to \cite{BZ05} for more details.

  Given a cluster algebra $\cA(\bx,\tilde B)$ of full rank, any choice of compatible Poisson structure gives rise to a canonical quantization of $\cA(\bx,\tilde B)$ by the following procedure, where $v$ denotes a formal variable:
  \begin{itemize}
    \item Each cluster should be replaced by a collection of quasi-commuting elements $\{X_1,\ldots,X_m\}$ satisfying $X_iX_j=v^{2\Lambda(\varepsilon_j,\varepsilon_i)}X_jX_i$ for a skew-symmetric bilinear form $\Lambda:\ZZ^m\times\ZZ^m\to\ZZ$, where $\varepsilon_i$ denotes the $i$th standard basis vector of $\ZZ^m$.  That is, this \emph{quantum cluster} generates a \emph{quantum torus} $\cT=\ZZ[v^{\pm1}]\langle X_1^{\pm1},\ldots,X_m^{\pm1}\rangle$ with a $\ZZ[v^{\pm1}]$-basis 
    \[\left\{X^\ba:=v^{\sum\limits_{i<j}a_ia_j\Lambda(\varepsilon_i,\varepsilon_j)}X_1^{a_1}\cdots X_m^{a_m}\Big|\ba=(a_1,\ldots,a_m)\in\ZZ^m\right\}.\]
    \item Each quantum cluster obtained by mutation from this one should again generate a quantum torus.  This forces a compatibility condition on the exchange matrix $\tilde B$ and the commutation matrix $\Lambda$ (by a slight abuse of notation we freely identify the bilinear form $\Lambda$ with a skew-symmetric matrix), that is 
    \begin{equation}\label{eq:quantum compatibility}
      \Lambda(\bb_k,\varepsilon_\ell)=0\quad\text{if $k\ne\ell$,}
    \end{equation}
    where $\bb_k$ denotes the $k^{th}$ column of $\tilde B$ thought of as an element of $\ZZ^m$.
    \item The above conditions still leave a great deal of freedom.  To restrict further, notice that the quasi-commutation relation above can be rewritten as $v^{\Lambda(\varepsilon_i,\varepsilon_j)}X_iX_j=v^{\Lambda(\varepsilon_j,\varepsilon_i)}X_jX_i$, in particular this implies that $\cT$ admits an anti-automorphism $X\mapsto \overline{X}$ which fixes each $X_i$ and sends $v$ to $v^{-1}$ (recall that $\Lambda$ was assumed skew-symmetric).  Since each initial cluster variable was \emph{bar-invariant}, the same should be true of all cluster variables since a cluster algebra is independent of the choice of initial cluster.  This uniquely determines how to incorporate powers of $v$ into the exchange relations \eqref{eq:exchange relation}; more precisely, the new variable obtained by mutation in direction $k$ must be
    \[X_k'=X^{-\varepsilon_k+\sum\limits_{b_{ik}>0}b_{ik}\varepsilon_i}+X^{-\varepsilon_k-\sum\limits_{b_{ik}<0}b_{ik}\varepsilon_i}.\]
  \end{itemize}
  \begin{exercise}
    Verify the compatibility condition \eqref{eq:quantum compatibility} and show that it is the same as the compatibility condition from Lemma~\ref{le:compatibility}.
  \end{exercise}

\section{Applications to Integrable Systems}\label{sec:integrable_systems}
	We now explore connections between cluster algebras and discrete integrable systems.  In this context, it is common to represent the exchange matrix as a quiver.

	\begin{definition}
		A \emph{quiver}\index{quiver} is a directed graph without loops (arrows from a vertex to itself) and without oriented $2$-cycles.
	\end{definition}
	
	A skew-symmetric $n\times n$ exchange matrix $B$ can be encoded as a quiver $Q$ on vertex set $\{1,2,\ldots, n\}$.  More precisely, if $i,j$ are such that $b_{ij}>0$ then $Q$ has $b_{ij}$ arrows from $i$ to $j$ (and none from $j$ to $i$).  Matrix mutation on $B$ induces so-called \emph{quiver mutation} on $Q$.
	
	\begin{definition}
		Given some $k=1,2,\ldots, n$, the \emph{quiver mutation} $\mu_k$ gives rise to $Q' = \mu_k(Q)$ by applying the following steps to $Q$:
		\begin{enumerate}
		\item for each pair of vertices $i,j$ and arrows $i\to k$ and $k \to j$, add an arrow $i \to j$,
		\item reverse all arrows $i \to k$ and $k \to j$,
		\item erase in turn pairs of arrows $i \to j$ and $j \to i$ to eliminate any $2$-cycles.
		\end{enumerate}
	\end{definition}
	
	Given a quiver $Q$ and a sequence $k_1,k_2,\ldots, k_l$ of its vertices (possibly with repeats), one can consider the dynamical system obtained by repeated applications of the composite mutation $\mu_{k_l} \circ \cdots \circ \mu_{k_2} \circ \mu_{k_1}$ starting from an initial seed $(\x, Q)$ or $Y$-seed $(\y,Q)$.  In the special case when the final quiver after one application equals the initial one, i.e.
	\begin{displaymath}
		\mu_{k_l} \circ \cdots \circ \mu_{k_2} \circ \mu_{k_1}(Q) = Q,
	\end{displaymath}
	the system amounts to iterating a fixed birational transformation.  It is natural to investigate these mappings for properties such as periodicity and integrability.  
	
	We focus on two (families of) quivers, one giving rise to the $Y$-system of type $A \times A$ and the other corresponding to the pentagram map. Some of the many other systems that have been investigated from a cluster algebra point of view are $Q$-systems \cite{dFK10}, $T$-systems \cite{dFK09}, Somos sequences \cite{H07}, and mutations in box products of quivers \cite{P15}. A systematic study of quivers that return to (isomorphic copies of) themselves under mutation was undertaken in \cite{FM11}.  We also refer the reader to \cite{N11,GNR16}, where periodicities of cluster algebra mutations are used to derive dilogarithm identities.

For the remainder of this section, all mutations are understood to be of $Y$-pattern type.

\subsection{Zamolodchikov periodicity}\index{Zamolodchikov periodicity}
	Fix positive integers $r$ and $s$ and consider the quiver $Q$ on vertex set $\{1,\ldots, r \} \times \{1,\ldots, s\}$ as illustrated in the case $r=4$, $s=3$ by
	
	\begin{pspicture}(5,4)
		\cnode(1,1){0.1}{v11}
		\cnode(2,1){0.1}{v21}
		\cnode(3,1){0.1}{v31}
		\cnode(4,1){0.1}{v41}
		\cnode(1,2){0.1}{v12}
		\cnode(2,2){0.1}{v22}
		\cnode(3,2){0.1}{v32}
		\cnode(4,2){0.1}{v42}
		\cnode(1,3){0.1}{v13}
		\cnode(2,3){0.1}{v23}
		\cnode(3,3){0.1}{v33}
		\cnode(4,3){0.1}{v43}
		\ncline{->}{v11}{v21}
		\ncline{->}{v31}{v21}
		\ncline{->}{v31}{v41}
		\ncline{<-}{v12}{v22}
		\ncline{<-}{v32}{v22}
		\ncline{<-}{v32}{v42}
		\ncline{->}{v13}{v23}
		\ncline{->}{v33}{v23}
		\ncline{->}{v33}{v43}
		\ncline{<-}{v11}{v12}
		\ncline{<-}{v13}{v12}
		\ncline{->}{v21}{v22}
		\ncline{->}{v23}{v22}
		\ncline{<-}{v31}{v32}
		\ncline{<-}{v33}{v32}
		\ncline{->}{v41}{v42}
		\ncline{->}{v43}{v42}
	\end{pspicture}
	
	Let $\mu_{\textrm{even}}$ be the compound mutation given by applying each $\mu_{i,j}$ with $i+j$ even in turn.  As no arrows connect these even vertices, it follows that these mutations commute and so the order does not matter.  Define $\mu_{\textrm{odd}}$ similarly as the compound mutation at all odd vertices.
	
	\begin{exercise}
		Let $Y = (Y_{i,j})_{i=1,\ldots, r; j=1,\ldots, s}$.  Then $\mu_{\textrm{even}}(Y,Q) = (Y',-Q)$ where
		\begin{equation} \label{eq:YSystemCases}
			Y_{i,j}' = \begin{cases}
			Y_{i,j}^{-1} & \textrm{if $i+j$ even} \\
			Y_{i,j}\frac{\displaystyle\prod_{|i-i'|=1}(1+Y_{i',j})}{\displaystyle\prod_{|j-j'|=1}(1+Y_{i,j'}^{-1})} & \textrm{if $i+j$ odd} \\
			\end{cases}
		\end{equation}
		and $-Q$ denotes $Q$ with all its arrows reversed.
	\end{exercise}
	
	Now begin with the $Y$-seed $(Y_0,Q)$ and let 
	\begin{align*}
	(Y_1,-Q) &= \mu_{\textrm{even}}(Y_0,Q) \\ 
	(Y_2,Q) &= \mu_{\textrm{odd}}(Y_1,-Q) \\ 
	(Y_3,-Q) &= \mu_{\textrm{even}}(Y_2,Q) \\ 
	&\vdots
	\end{align*}
	where $Y_t = (Y_{ijt})_{i=1,\ldots, r; j=1,\ldots, s}$.  It follows from the above that
	\begin{displaymath}
		Y_{i,j,t+1} = \begin{cases}
		Y_{i,j,t}^{-1} & \textrm{if $i+j+t$ even} \\
		Y_{i,j,t}\frac{\displaystyle\prod_{|i-i'|=1}(1+Y_{i',j,t})}{\displaystyle\prod_{|j-j'|=1}(1+Y_{i,j',t}^{-1})} & \textrm{if $i+j+t$ odd} \\
		\end{cases}
	\end{displaymath}
	It is convenient to consider only the $Y_{i,j,t}$ with $i+j+t$ even, which satisfy their own recurrence
	\begin{equation} \label{eq:YSystem}
		Y_{i,j,t-1}Y_{i,j,t+1} = \frac{\displaystyle\prod_{|i-i'|=1}(1+Y_{i',j,t})}{\displaystyle\prod_{|j-j'|=1}(1+Y_{i,j',t}^{-1})}
	\end{equation}
	for all $i+j+t$ odd.  This recurrence is called the type $A_r \times A_s$ $Y$-system, a special case of a family of systems conjectured by Zamolodchikov to be periodic.  The proof, in this case, is due to A. Volkov \cite{V07}.
	
	\begin{theorem}
		The $Y$-system \eqref{eq:YSystem} on variables $Y_{i,j,t}$ with $i+j+t$ odd, $i=1,\ldots, r$, and $j=1,\ldots, s$ has period $2(r+s+2)$ in the $t$-direction, that is
		\begin{displaymath}
			Y_{i,j,t+2(r+s+2)} = Y_{i,j,t}.
		\end{displaymath}
	\end{theorem}
	
	Returning to the $Y$-pattern point of view, the initial seed $(Y_0,Q)$ consists of the $Y_{i,j,0}$ for $i+j$ even and the $Y_{i,j,-1}^{-1}$ for $i+j$ odd.  The seed $(Y_2,Q) = \mu_{\textrm{odd}}(\mu_{\textrm{even}}(Y_0,Q))$ consists of the $Y_{i,j,2}$ for $i+j$ even and the $Y_{i,j,1}^{-1}$ for $i+j$ odd.  The periodicity theorem, then, asserts that the rational map $\mu_{\textrm{odd}} \circ \mu_{\textrm{even}}$ has order $r+s+2$.
	
	\subsection{The pentagram map} 
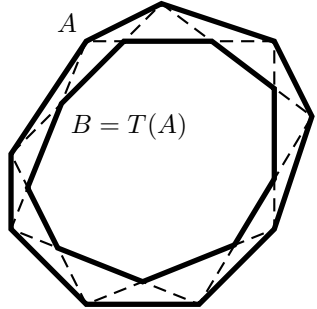
\begin{figure} \label{fig:pentagram}
	\begin{pspicture}(5,5)
	\pspolygon[linewidth=2pt](1,2)(1,3)(2,4.5)(3,5)(4.5,4.5)(5,3.5)(4.5,2)(3.5,1)(2,1)
  \pspolygon[linestyle=dashed](1,2)(2,4.5)(4.5,4.5)(4.5,2)(2,1)(1,3)(3,5)(5,3.5)(3.5,1)
  \pspolygon[linewidth=2pt](1.22,2.55)(1.67,3.67)(2.5,4.5)(3.67,4.5)(4.5,3.87)(4.5,2.67)(3.97,1.79)(2.75,1.3)(1.62,1.75)
	\uput[ul](2,4.5){$A$}
	\uput[dr](1.67,3.67){$B=T(A)$}
	\end{pspicture}
	\caption{An application of the pentagram map}
	\end{figure}

The pentagram map\index{pentagram map} is a discrete dynamical system defined on the space of polygons in the projective plane.  Figure \ref{fig:pentagram} shows a polygon $A$ and the corresponding output $B=T(A)$, with $T$ being the notation for the pentagram map.  Each vertex of $B$ lies at the intersection of two consecutive ``shortest'' diagonals of $A$.
	
	The pentagram map was introduced by R. Schwartz \cite{S92} and extensively studied by V. Ovsienko, Schwartz and S. Tabachnikov \cite{OST10, OST13}.  In particular, they demonstrated that the pentagram map is a completely integrable system.  Another take on integrability was provided by F. Soloviev \cite{S13}.  In the same span of time, M. Glick \cite{G11} described the pentagram map as mutations in a $Y$-pattern, building on work of Schwartz \cite{S08} who had found a lift to the octahedron recurrence.  Finally, M. Gekhtman, M. Shapiro, Tabachnikov and A. Vainshtein \cite{GSTV12} gave a uniform treatment of integrability and the cluster algebra structure of the pentagram map, and generalizations thereof, in terms of weighted networks on tori.  A similar framework is provided by the cluster integrable systems of A. Goncharov and R. Kenyon \cite{GK13}, albeit without any explicit mention of the pentagram map.
	
	The connection between the pentagram map and $Y$-patterns comes by way of certain geometrically defined coordinates on the space of polygons.  The \emph{cross ratio}\index{cross ratio} of $4$ real numbers $x_1,x_2,x_3,x_4$ is 
	\begin{displaymath}
		\chi(x_1,x_2,x_3,x_4) = \frac{(x_1-x_2)(x_3-x_4)}{(x_1-x_3)(x_2-x_4)}.
	\end{displaymath}
	The cross ratio is invariant under projective transformations, that is 
	\begin{displaymath}
		\chi(f(x_1),f(x_2),f(x_3),f(x_4)) = \chi(x_1,x_2,x_3,x_4)
	\end{displaymath}
	for any fractional linear map $f(x) = (ax+b)/(cx+d)$.  For this reason, there is a well-defined notion in the plane of the cross ratio of 4 collinear points or of 4 concurrent lines.
	
\begin{figure} 
\begin{pspicture}(6,4)
	\rput(0,-0.5){
  \pspolygon[showpoints=true,linewidth=2pt](1,1)(4,1)(5,3)(3,4)(1,3)
	\uput[l](1,2){$1$}
  \uput[225](1,1){$2$}
	\uput[d](2.5,1){$3$}
  \uput[270](4,1){$4$}
	\uput[dr](4.5,2){$5$}
  \uput[0](5,3){$6$}
	\uput[ur](4,3.5){$7$}
  \uput[90](3,4){$8$}
	\uput[ul](2,3.5){$9$}
  \uput[l](1,3){$10$}
  \pspolygon[linestyle=dashed](1,1)(5,3)(1,3)(4,1)(3,4)
  \pspolygon[showpoints=true,linewidth=2pt](1.92,2.38)(2.71,1.86)(3.57,2.28)(3.33,3)(2.33,3)
  \uput[l](1.92,2.38){$1$}
	\uput[dl](2.31,2.12){$2$}
	\uput[d](2.71,1.86){$3$}
	\uput[dr](3.14,2.07){$4$}
  \uput[330](3.57,2.28){$5$}
	\uput[20](3.45,2.64){$6$}
  \uput[ur](3.33,3){$7$}
	\uput[u](2.83,3){$8$}
  \uput[ul](2.33,3){$9$}
	\uput[160](2.12,2.69){$10$}
  
  }
\end{pspicture}
\caption{Possible labelings of two polygons related by the pentagram map}
\label{fig:labeling}
\end{figure}
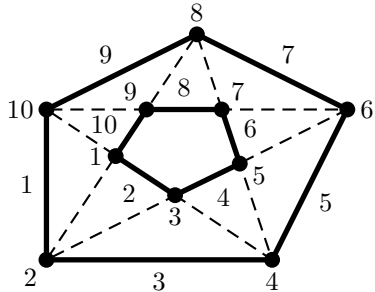

	Given an $n$-gon $A$, label its sides and vertices consecutively with the integers $\{1,2,\ldots, 2n\}$.  Such a labeling induces a canonical labeling on $T(A)$ as is illustrated in Figure \ref{fig:labeling}.  Note that the parities of the vertex labels and of the edge labels are interchanged by $T$.  In what follows, let $\join{i}{j}$ denote the line between two points labeled $i$ and $j$ and let $\meet{k}{l}$ denote the point of intersection of two lines labeled $k$ and $l$.
	
\begin{figure}
\psset{unit=.8cm}
\begin{pspicture}(12,6)
\rput(0,.5){
\psline(1,1)(3,1)(5,2)(5,3)(4,4)
\uput[d](1,1){$1$}
\uput[d](2,1){$2$}
\uput[d](3,1){$3$}
\uput[dr](4,1.5){$4$}
\uput[dr](5,2){$5$}
\uput[r](5,2.5){$6$}
\uput[r](5,3){$7$}
\uput[ur](4.5,3.5){$8$}
\uput[ur](4,4){$9$}
\psline[linewidth=2pt](1,1)(5,2)
\psline[linewidth=2pt](3,1)(5,2)
\psline[linewidth=2pt](5,3)(5,2)
\psline[linewidth=2pt](4,4)(5,2)
}
\rput(7,0){
\psline(1,1)(3,1)(5,2)(5,3)(4,4)(2,5)
\uput[d](1,1){$1$}
\uput[d](2,1){$2$}
\uput[d](3,1){$3$}
\uput[345](4,1.5){$4$}
\uput[r](5,2){$5$}
\uput[r](5,2.5){$6$}
\uput[r](5,3){$7$}
\uput[dl](4.5,3.5){$8$}
\uput[u](4,4){$9$}
\uput[u](3,4.5){$10$}
\uput[u](2,5){$11$}
\psline(3,1)(5,1)(5,3.5)(4,4)
\psdots[dotsize=2pt 3](5,1)(5,2)(5,3)(5,3.5)
}
\end{pspicture}
\psset{unit=1cm}
\caption{The cross ratios corresponding to the two types of $y$-parameters.  On the left, $-(y_5)^{-1}$ is the cross ratio of the lines $\join{5}{1}$, $4$, $6$, and $\join{5}{9}$.  On the right, $-y_6$ is the cross ratio of the points $\meet{6}{2}$, $5$, $7$, and $\meet{6}{10}$ .} 
\label{fig:defy}
\end{figure}
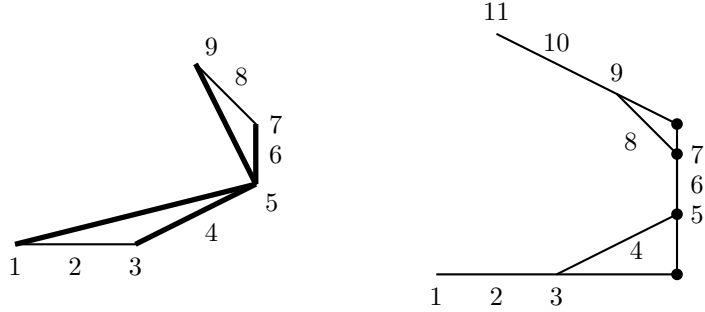

\begin{definition}
	The \emph{$y$-parameters} of a polygon $A$ are real numbers $y_1(A),\ldots, y_{2n}(A)$ defined by
	\begin{align*}
		y_i(A) = \begin{cases}
		-\left(\chi(\join{i}{(i-4)}, i-1, i+1, \join{i}{(i+4)})\right)^{-1} & \textrm{ if $i$ is a vertex of $A$} \\
		-\chi(\meet{i}{(i-4)}, i-1, i+1, \meet{i}{(i+4)}) & \textrm{ if $i$ is a side of $A$} \\
		\end{cases}
	\end{align*}
	This definition is illustrated in Figure \ref{fig:defy}.
\end{definition}

	\begin{proposition}
		Let $A$ be an $n$-gon with $y$-parameters $y_i = y_i(A)$.  Let $B= T(A)$ and denote its $y$-parameters $y_i' = y_i(B)$.  Then
		\begin{equation} \label{eq:T}
			y_i' = \begin{cases}
			y_i^{-1} & \textrm{if $i$ is a side of $B$} \\
			y_i\frac{(1+y_{i-1})(1+y_{i+1})}{(1+y_{i-3}^{-1})(1+y_{i+3}^{-1})} & \textrm{if $i$ is a vertex of $B$}
			\end{cases}
		\end{equation}
		In this equation, all indices are considered modulo $2n$.
	\end{proposition}
	
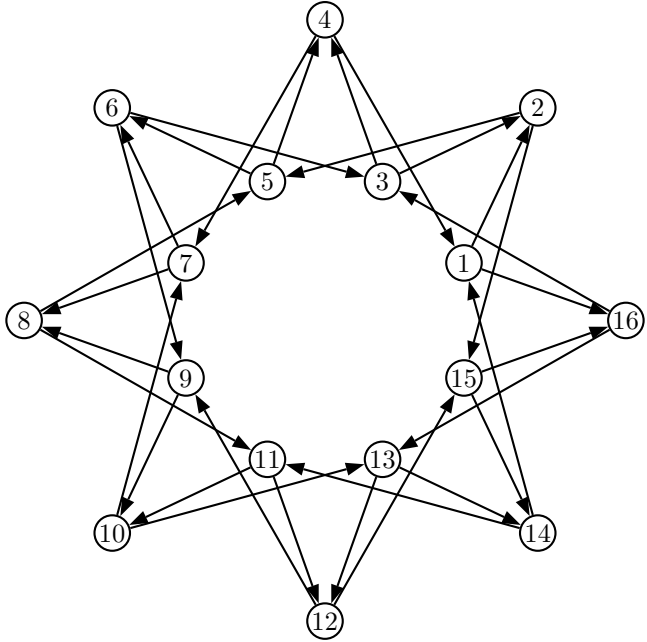
\begin{figure}
\begin{pspicture}(-5,-5)(5,5)
\SpecialCoor
\cnode(2;22.5){.25}{v1}
\rput(v1){1}
\cnode(4;45){.25}{v2}
\rput(v2){2}
\cnode(2;67.5){.25}{v3}
\rput(v3){3}
\cnode(4;90){.25}{v4}
\rput(v4){4}
\cnode(2;112.5){.25}{v5}
\rput(v5){5}
\cnode(4;135){.25}{v6}
\rput(v6){6}
\cnode(2;157.5){.25}{v7}
\rput(v7){7}
\cnode(4;180){.25}{v8}
\rput(v8){8}
\cnode(2;202.5){.25}{v9}
\rput(v9){9}
\cnode(4;225){.25}{v10}
\rput(v10){10}
\cnode(2;247.5){.25}{v11}
\rput(v11){11}
\cnode(4;270){.25}{v12}
\rput(v12){12}
\cnode(2;292.5){.25}{v13}
\rput(v13){13}
\cnode(4;315){.25}{v14}
\rput(v14){14}
\cnode(2;337.5){.25}{v15}
\rput(v15){15}
\cnode(4;0){.25}{v16}
\rput(v16){16}

\psset{arrowsize=5pt}
\psset{arrowinset=0}
\ncline{->}{v1}{v16}   \ncline{->}{v1}{v2}
\ncline{->}{v3}{v2}    \ncline{->}{v3}{v4} 
\ncline{->}{v5}{v4}    \ncline{->}{v5}{v6}
\ncline{->}{v7}{v6}    \ncline{->}{v7}{v8}
\ncline{->}{v9}{v8}    \ncline{->}{v9}{v10}
\ncline{->}{v11}{v10}  \ncline{->}{v11}{v12}
\ncline{->}{v13}{v12}  \ncline{->}{v13}{v14}
\ncline{->}{v15}{v14}  \ncline{->}{v15}{v16}
\ncline{->}{v16}{v13} \ncline{->}{v16}{v3}
\ncline{->}{v2}{v15} \ncline{->}{v2}{v5}
\ncline{->}{v4}{v1} \ncline{->}{v4}{v7}
\ncline{->}{v6}{v3} \ncline{->}{v6}{v9}
\ncline{->}{v8}{v5} \ncline{->}{v8}{v11}
\ncline{->}{v10}{v7} \ncline{->}{v10}{v13}
\ncline{->}{v12}{v9} \ncline{->}{v12}{v15}
\ncline{->}{v14}{v11} \ncline{->}{v14}{v1}

\end{pspicture}
\caption{The quiver related to the action of the pentagram map on octagons} \label{fig:GlickQuiver}
\end{figure}

	The formula \eqref{eq:T} bears a distinct resemblance to \eqref{eq:YSystemCases}, so it is plausible that it can also be described by $Y$-pattern mutations for an appropriate quiver.  The desired quiver $Q_n$ has vertex set $\{1,2,\ldots,2n\}$ and arrows $j \to (j\pm1)$ and $j \leftarrow (j \pm 3)$ for each odd $j$, with vertices considered modulo $2n$.  The quiver $Q_8$ is shown in Figure~\ref{fig:GlickQuiver}.  Similar to before, define compound mutations
	\begin{align*}
		\mu_{\textrm{odd}} &= \mu_{2n-1} \circ \cdots \circ \mu_3 \circ \mu_1, \\
		\mu_{\textrm{even}} &= \mu_{2n} \circ \cdots \circ \mu_4 \circ \mu_2. 
	\end{align*}
	
\begin{theorem}[\cite{G11}] \label{thm:pentagram}
	Let $A$ be an $n$-gon with $y$-parameters $y_1,\ldots, y_{2n}$ and let $k$ be a positive integer.  Beginning from the $Y$-seed $((y_1,\ldots, y_{2n}), Q_n)$, apply $k$ compound mutations alternating between $\mu_{\textrm{even}}$ and $\mu_{\textrm{odd}}$.  The result will have the form $((y_1',\ldots, y_{2n}'), (-1)^kQ_n)$ where 
	\begin{displaymath}
	y_i' = y_i(T^k(A)). 
	\end{displaymath}
\end{theorem}

\begin{remark}
	It is natural to ask what values the $n$-tuple $(y_1(A),\ldots, y_{2n}(A))$ can take, and also to what extent this data suffices to reconstruct $A$.  These questions are easier to answer after extending to a larger family of objects called twisted polygons.  In this setting, the $y_i$ satisfy a single relation $y_1\cdots y_{2n}=1$.  Moreover, the $y_i$ determine a twisted polygon uniquely up to projective equivalence and a one-parameter rescaling operation due to Schwartz \cite{S08}.  The $Y$-pattern dynamics described in Theorem \ref{thm:pentagram}, then, characterize the dynamics of the pentagram map on the space of twisted polygons modulo these equivalences.
\end{remark}
	
	We now give a somewhat informal presentation of integrability of the pentagram map from the cluster algebra perspective.  The first ingredient is a compatible Poisson structure, which can be obtained directly from the $Y$-pattern analogue of Theorem \ref{thm:Poisson}.
	
	\begin{proposition}[\cite{OST10, GSTV12}] \label{prop:pentagramBracket}
		Define a Poisson bracket by $\{y_i,y_j\} = b_{ij}y_iy_j$ where $B$ is the exchange matrix associated to the pentagram quiver $Q_n$.  Then this bracket is preserved by the pentagram map, i.e.
		\begin{displaymath}
			\{f,g\} \circ T = \{f\circ T, g \circ T\}
		\end{displaymath}
		for all functions $f$ and $g$ of the $y$-parameters.
	\end{proposition}
	
\begin{figure}
	\begin{pspicture}(0,.5)(10,3)
	\rput(1,1){\rnode{v1}{$1$}}
	\rput(2,1){\rnode{v2}{$2$}}
	\rput(3,1){\rnode{v3}{$3$}}
	\rput(4,1){\rnode{v4}{$4$}}
	\rput(5,1){\rnode{v5}{$5$}}
	\rput(6,1){\rnode{v6}{$6$}}
	\rput(7,1){\rnode{v7}{$7$}}
	\rput(8,1){\rnode{v8}{$8$}}
	\rput(9,1){\rnode{v9}{$1$}}
	\rput(1,2){\rnode{w1}{$4$}}
	\rput(2,2){\rnode{w2}{$5$}}
	\rput(3,2){\rnode{w3}{$6$}}
	\rput(4,2){\rnode{w4}{$7$}}
	\rput(5,2){\rnode{w5}{$8$}}
	\rput(6,2){\rnode{w6}{$1$}}
	\rput(7,2){\rnode{w7}{$2$}}
	\rput(8,2){\rnode{w8}{$3$}}
	\rput(9,2){\rnode{w9}{$4$}}
	\rput(1,3){\rnode{x1}{$7$}}
	\rput(2,3){\rnode{x2}{$8$}}
	\rput(3,3){\rnode{x3}{$1$}}
	\rput(4,3){\rnode{x4}{$2$}}
	\rput(5,3){\rnode{x5}{$3$}}
	\rput(6,3){\rnode{x6}{$4$}}
	\rput(7,3){\rnode{x7}{$5$}}
	\rput(8,3){\rnode{x8}{$6$}}
	\rput(9,3){\rnode{x9}{$7$}}
	\ncline{->}{v1}{v2}
	\ncline{->}{v3}{v2}
	\ncline{->}{v3}{v4}
	\ncline{->}{v5}{v4}
	\ncline{->}{v5}{v6}
	\ncline{->}{v7}{v6}
	\ncline{->}{v7}{v8}
	\ncline{->}{v9}{v8}
	\ncline{<-}{w1}{w2}
	\ncline{<-}{w3}{w2}
	\ncline{<-}{w3}{w4}
	\ncline{<-}{w5}{w4}
	\ncline{<-}{w5}{w6}
	\ncline{<-}{w7}{w6}
	\ncline{<-}{w7}{w8}
	\ncline{<-}{w9}{w8}
	\ncline{->}{x1}{x2}
	\ncline{->}{x3}{x2}
	\ncline{->}{x3}{x4}
	\ncline{->}{x5}{x4}
	\ncline{->}{x5}{x6}
	\ncline{->}{x7}{x6}
	\ncline{->}{x7}{x8}
	\ncline{->}{x9}{x8}
	\ncline{<-}{w1}{w2}
	\ncline{->}{w1}{v1}
	\ncline{<-}{w2}{v2}
	\ncline{->}{w3}{v3}
	\ncline{<-}{w4}{v4}
	\ncline{->}{w5}{v5}
	\ncline{<-}{w6}{v6}
	\ncline{->}{w7}{v7}
	\ncline{<-}{w8}{v8}
	\ncline{->}{w9}{v9}
	\ncline{->}{w1}{x1}
	\ncline{<-}{w2}{x2}
	\ncline{->}{w3}{x3}
	\ncline{<-}{w4}{x4}
	\ncline{->}{w5}{x5}
	\ncline{<-}{w6}{x6}
	\ncline{->}{w7}{x7}
	\ncline{<-}{w8}{x8}
	\ncline{->}{w9}{x9}

\end{pspicture}
	\caption{A lift of $Q_4$ to a doubly periodic quiver in the plane.}
	\label{fig:quiverLift}
	\end{figure}
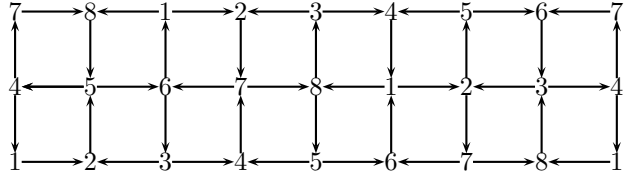

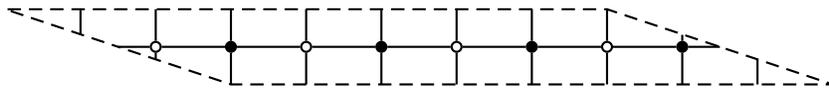
\begin{figure}
	\begin{pspicture}(-1,0)(10,2)
\cnode(1,1){.08}{v1}
\cnode*(2,1){.08}{v2}
\cnode(3,1){.08}{v3}
\cnode*(4,1){.08}{v4}
\cnode(5,1){.08}{v5}
\cnode*(6,1){.08}{v6}
\cnode(7,1){.08}{v7}
\cnode*(8,1){.08}{v8}
\pnode(.5,1){w1}
\pnode(8.5,1){e8}
\pnode(0,1.5){n0}
\pnode(1,1.5){n1}
\pnode(2,1.5){n2}
\pnode(3,1.5){n3}
\pnode(4,1.5){n4}
\pnode(5,1.5){n5}
\pnode(6,1.5){n6}
\pnode(7,1.5){n7}
\pnode(8,1.16){n8}
\pnode(9,.83){n9}
\pnode(0,1.16){s0}
\pnode(1,.83){s1}
\pnode(2,.5){s2}
\pnode(3,.5){s3}
\pnode(4,.5){s4}
\pnode(5,.5){s5}
\pnode(6,.5){s6}
\pnode(7,.5){s7}
\pnode(8,.5){s8}
\pnode(9,.5){s9}
\pspolygon[linestyle=dashed](-1,1.5)(7,1.5)(10,.5)(2,.5)
\ncline{w1}{v1}
\ncline{v1}{v2}
\ncline{v2}{v3}
\ncline{v3}{v4}
\ncline{v4}{v5}
\ncline{v5}{v6}
\ncline{v6}{v7}
\ncline{v7}{v8}
\ncline{v8}{e8}
\ncline{v1}{n1}
\ncline{v2}{n2}
\ncline{v3}{n3}
\ncline{v4}{n4}
\ncline{v5}{n5}
\ncline{v6}{n6}
\ncline{v7}{n7}
\ncline{v8}{n8}
\ncline{v1}{s1}
\ncline{v2}{s2}
\ncline{v3}{s3}
\ncline{v4}{s4}
\ncline{v5}{s5}
\ncline{v6}{s6}
\ncline{v7}{s7}
\ncline{v8}{s8}
\ncline{s0}{n0}
\ncline{s9}{n9}
\end{pspicture}
	\caption{A bipartite graph on the torus dual to the quiver $Q_4$.}
	\label{fig:torus}
	\end{figure}

	Some more work is needed to get at the conserved quantities.  We follow the approach of \cite{GK13}.  The quiver $Q_n$ can be lifted to an infinite doubly periodic quiver in the plane as in Figure \ref{fig:quiverLift}.  This lift represents an embedding of $Q_n$ on a torus.  Let $G = (V,E)$ be the dual graph on the torus with $V = \{1,2,\ldots, 2n\}$ where $i \in V$ corresponds to the face of the quiver with vertices $i \pm 2, i \pm 1$.  Coincidentally, $G$ is identical to $Q_n$ except without orientations, so $\overline{ij} \in E$ if and only if $i-j \equiv \pm 1, \pm 3 \pmod{2n}$.  Figure \ref{fig:torus} shows the graph $G$ on the torus in the case $n=4$, with the vertices $1,2,\ldots, 8$ appearing in order from left to right.
	
	Assign edge weights to $G$ as follows: each ``horizontal'' edge $\overline{i(i+1)}$ gets weight $1$ and each ``vertical edge'' $\overline{(i-2)(i+1)}$ gets weight $(-1)^ix_i$.  The $x_i$ represent the \emph{corner invariants} \cite{S08} of a polygon, certain coordinates related to the $y$-parameters by
	\begin{displaymath}
		y_i = \begin{cases}
		-(x_ix_{i+1})^{-1} & \textrm{if $i$ is a vertex of $A$;} \\
		-x_ix_{i+1} & \textrm{if $i$ is a side of $A$.} \\
		\end{cases}
	\end{displaymath}

	A \emph{perfect matching} of $G$ is a collection $M$ of its edges such that each $i \in V$ is an endpoint of exactly one $e \in M$.  The \emph{weight} of $M$, denoted $wt(M)$, is the product of the weights of its edges.  For instance, 
	\begin{displaymath}
		M = \{\overline{12}, \overline{36}, \overline{47}, \overline{58}\}
	\end{displaymath}
	is a perfect matching of $G_4$ and its weight is $(1)(-x_5)(x_6)(-x_7) = x_5x_6x_7$.  The conserved quantities of the pentagram map are sums of weights of perfect matchings, refined by a notion of homology class of a matching defined in \cite{GK13}.
		
	\begin{proposition} \label{prop:pentagramIntegrals}
		For any homology class, the sum
		\begin{displaymath}
			\sum_M wt(M)
		\end{displaymath}
		over perfect matchings $M$ of $G_n$ in that class is a conserved quantity of the pentagram map acting on $n$-gons.	\end{proposition}
	
	In the case of the graph $G_4$ pictured in Figure \ref{fig:torus}, there are $22$ perfect matchings which represent $8$ different homology classes.  Table \ref{table:integrals} gives the $8$ corresponding conserved quantities for the pentagram map.  The first $2$ are trivial but the remaining $6$ are (up to a sign) the original conserved quantities discovered by Schwartz \cite{S08} by different means.
	
	\begin{table}[ht]
	\begin{tabular}{l|l|l}
	Name & Formula & Representative matching\\
	\hline
	n.a.  & $1$ & $\{\overline{12}, \overline{34}, \overline{56}, \overline{78}\}$\\
	n.a.  & $1$ & $\{\overline{23}, \overline{45}, \overline{67}, \overline{18}\}$\\
	$O_1$ & $-x_1-x_3-x_5-x_7+x_1x_2x_3+x_3x_4x_5$ & $\{\overline{27}, \overline{34}, \overline{56}, \overline{18}\}$ \\
	& $+x_5x_6x_7+x_1x_7x_8$ \\
	$E_1$ & $x_2+x_4+x_6+x_8-x_2x_3x_4-x_4x_5x_6$ & $\{\overline{38}, \overline{12}, \overline{45}, \overline{67}\}$ \\
	& $-x_6x_7x_8-x_1x_2x_8$ \\
	$O_2$ & $x_1x_5 + x_3x_7$ & $\{\overline{27}, \overline{36}, \overline{45}, \overline{18}\}$ \\
	$E_2$ & $x_2x_6 + x_4x_8$ & $\{\overline{38}, \overline{47}, \overline{12}, \overline{56}\}$ \\
	$O_4$ & $x_1x_3x_5x_7$ & $\{\overline{14}, \overline{36}, \overline{58}, \overline{27}\}$ \\
	$E_4$ & $x_2x_4x_6x_8$ & $\{\overline{16}, \overline{38}, \overline{25}, \overline{47}\}$ \\
	\end{tabular}
	\caption{Conserved quantities for the pentagram map acting on twisted quadrilaterals.  The matching $M$ is such that $wt(M)$ is the first term appearing in the corresponding formula.}
	\label{table:integrals}
	\end{table}
	
	The invariant Poisson bracket from Proposition \ref{prop:pentagramBracket} and the conserved quantities from Proposition \ref{prop:pentagramIntegrals} are the main ingredients for complete integrability, but there are several more conditions that they need to satisfy.  We conclude by listing these properties, proofs for which can be found in \cite{OST10}, \cite{GSTV12} and \cite{GK13}.
	\begin{itemize}
	\item The conserved quantities are algebraically independent.
	\item The conserved quantities Poisson commute, i.e. the bracket of any two of them equals zero.
	\item The number of Casimirs plus twice the number of other conserved quantities equals the dimension of the system.
	\end{itemize}
	In the last item, the term \emph{Casimir}\index{Casimir} refers to a function that Poisson commutes with all other functions.  In the case summarized in Table \ref{table:integrals}, the $4$ functions $O_2,E_2,O_4,E_4$ are Casimirs while the other $2$ conserved quantities $O_1$ and $E_1$ are not.  The calculation checks out as $4+2(2)=8$ is the dimension of the space of twisted quadrilaterals.
	
	\medskip
	\textbf{Acknowledgments.} These notes are based on lectures we presented during the 2016 ASIDE summer school in Montreal.  We thank the ASIDE organizers for inviting us and also for coordinating the efforts to compile the lecture notes.  We thank Sophie Morier-Genoud and the anonymous referee for several helpful comments.  The final version of the notes have been published in \cite{GR17}.
	



\begin{thebibliography}{99}
  	\bibitem{ABCGPT14} N. Arkani-Hamed, J. Bourjaily, F. Cachazo, A. Goncharov, A. Postnikov, and J. Trnka, Scattering amplitudes and the positive Grassmannian, arXiv preprint: 1212.5605v2.
  	\bibitem{BFZ05} A. Berenstein, S. Fomin, and A. Zelevinsky, Cluster algebras III: upper bounds and double Bruhat cells, \textsl{J. London Math. Soc.} \textbf{14} (1976), 183–-187.
  	\bibitem{BR15} A. Berenstein and D. Rupel, Quantum cluster characters of Hall algebras, \textsl{Sel. Math. New Ser.} \textbf{21} (2015), no. 4, 1121--1176.
  	\bibitem{BZ05} A. Berenstein and A. Zelevinsky, Quantum cluster algebras. \textsl{Adv. Math.} \textbf{195} (2005), no. 2, 405--455.
	\bibitem{BMRRT06} A. Buan, R. Marsh, I. Reiten, M. Reineke, and G. Todorov, Tilting theory and cluster combinatorics, \textsl{Adv. Math.} \textbf{204} (2006), no. 2, 572–-618.
  	\bibitem{CC06} P. Caldero and F. Chapoton, Cluster algebras as Hall algebras of quiver representations, \textsl{Comment. Math. Helv.} \textbf{81} (2006), no. 3, 595--616.
  	\bibitem{CK06} P. Caldero and B. Keller, From triangulated categories to cluster algebras II, \textsl{Ann. Sci. \'Ecole Norm. Sup.} (4) \textbf{39} (2006), no. 6, 983--1009.
	\bibitem{dFK09} P. DiFrancesco and R. Kedem, Positivity of the $T$-system cluster algebra, \textsl{Electron. J. Combin.} \textbf{16} (2009), 39 pp.
	\bibitem{dFK10} P. DiFrancesco and R. Kedem, $Q$-systems, heaps, paths and cluster positivity, \textsl{Comm. Math. Phys.} \textbf{293} (2010), 727--802.
  	\bibitem{EF12} R. Eager and S. Franco, Colored BPS pyramid partition functions, quivers and cluster transformations, \textsl{J. High Energy Phys.} (2012), no. 9, 038.
	\bibitem{FG06} V. Fock and A. Goncharov, Moduli spaces of local systems and higher Teichm\"uller theory, \textsl{Publ. Math. Inst. Hautes \'Etudes Sci.} \textbf{103} (2006), 1--211.
	\bibitem{FST08} S. Fomin, M. Shapiro, and D. Thurston, Cluster algebras and triangulated surfaces. I. Cluster complexes, \textsl{Acta Math.} \textbf{201} (2008), 83--146.
	\bibitem{FZ99} S. Fomin and A. Zelevinsky, Double Bruhat cells and total positivity, \textsl{J. Amer. Math. Soc.} \textbf{12} (1999), 335--380.
	\bibitem{FZ02} S. Fomin and A. Zelevinsky, Cluster algebras I: Foundations, \textsl{J. Amer. Math. Soc.} \textbf{15} (2002), 497--529.
	\bibitem{FZ02b} S. Fomin and A. Zelevinsky, The Laurent phenomenon, \textsl{Adv. in Applied Math.} \textbf{28} (2002), 119-–144.
  	\bibitem{FZ03} S. Fomin and A. Zelevinsky, Cluster algebras II: Finite type classification, \textsl{Invent. Math.} \textbf{154} (2003), 63-–121.
	\bibitem{FZ07} S. Fomin and A. Zelevinsky, Cluster algebras IV: Coefficients, \textsl{Compos. Math.} \textbf{143} (2007), 112--164.
	\bibitem{FM11} A. Fordy and R. Marsh, Cluster mutation-periodic quivers and associated Laurent sequences, \textsl{J. Algebraic Combin.} \textbf{34} (2011), no. 1, 19--66.
	\bibitem{GLS06} C. Geiss, B. Leclerc, and J. Schr\"oer, Rigid modules over preprojective algebras, \textsl{Invent. Math} \textbf{165} (2006), 589--632.
 	\bibitem{GLS13} C. Geiss, B. Leclerc, and J. Schr\"oer, Cluster structures on quantum coordinate rings, \textsl{Selecta Math. (N.S.)} \textbf{19} (2013), no. 2, 337--397.
  	\bibitem{GLS} C. Geiss, B. Leclerc, and J. Schr\"oer, Cluster algebras in algebraic Lie theory, \textsl{Transf. Groups} \textbf{18} (2013), 149--178.
  	\bibitem{GNR16} M. Gekhtman, T. Nakanishi, D. Rupel, Hamiltonian and Lagrangian formalisms of mutations in cluster algebras and application to dilogarithm identities, arXiv preprint: 1611.02813.
	\bibitem{GSTV12} M. Gekhtman, M. Shapiro, S. Tabachnikov, and A. Vainshtein, Higher pentagram maps, weighted directed networks, and cluster dynamics, \textsl{Electron. Res. Announc. Math. Sci.} \textbf{19} (2012), 1--17.
	\bibitem{GSV03} M. Gekhtman, M. Shapiro, and A. Vainshtein,  Cluster algebras and Poisson geometry, \textsl{Mosc. Math. J.} \textbf{3} (2003), 899--934.
  	\bibitem{GSV10} M. Gekhtman, M. Shapiro, and A. Vainshtein,  Cluster algebras and Poisson geometry, Mathematical Surveys and Monographs, \textbf{167}. American Mathematical Society, Providence, RI, 2010.
	\bibitem{G11} M. Glick, The pentagram map and Y-patterns, \textsl{Adv. Math.} \textbf{227} (2011), 1019--1045.
	\bibitem{GR17} M. Glick and D. Rupel: Introduction to Cluster Algebras. In: Levi D., Rebelo R., Winternitz P. (eds) Symmetries and Integrability of Difference Equations. CRM Series in Mathematical Physics. Springer, Cham (2017).
	\bibitem{GK13} A. Goncharov and R. Kenyon, Dimers and cluster integrable systems, \textsl{Ann. Sci. \'Ec. Norm. Sup\'er.} \textbf{46} (2013), 747--813.
  	\bibitem{GHKK14} M. Gross, P. Hacking, S. Keel, and M. Kontsevich, Canonical bases for cluster algebras, arXiv preprint:1411.1394.
	\bibitem{H91} A. Hatcher, On triangulations of surfaces, \textsl{Topology Appl.} \textbf{40} (1991), 189--194.
	\bibitem{H07} A. Hone, Sigma function solution of the initial value problem for Somos 5 sequences, \textsl{Trans. Amer. Math. Soc.} \textbf{359} (2007), no. 10, 5019--5034.
  	\bibitem{K10} B. Keller, Categorification of acyclic cluster algebras: an introduction, arXiv preprint: 0801.3103v2.
  	\bibitem{K12} Y. Kimura, Quantum unipotent subgroup and dual canonical basis, \textsl{Kyoto J. Math.} \textbf{52} (2012), no. 2, 277--331.
  	\bibitem{KQ14} Y. Kimura and F. Qin, Graded quiver varieties, quantum cluster algebras and dual canonical basis, \textsl{Adv. Math.} \textbf{262} (2014), 261--312.
  	\bibitem{KZ02} M. Kogan and A. Zelevinsky, On symplectic leaves and integrable systems in standard complex semisimple Poisson-Lie groups, \textsl{Int. Math. Res. Not.} (2002), no. 32, 1685--1702.
  	\bibitem{LLZ14} K. Lee, L. Li, and A. Zelevinsky, Greedy elements in rank 2 cluster algebras, \textsl{Sel. Math. New Ser.} \textbf{20} (2014), 57--82.
  	\bibitem{LS13} K. Lee and R. Schiffler, A Combinatorial Formula for Rank 2 Cluster Variables, \textsl{J. Alg. Comb} \textbf{37} (2013), 67--85.
  	\bibitem{LS15} K. Lee and R. Schiffler, Positivity for cluster algebras, \textsl{Ann. Math.} \textbf{182} (2015), 73--125.
	\bibitem{MP07} G. Musiker and J. Propp, Combinatorial interpretations for rank-two cluster algebras of affine type, \textsl{ Electron. J. Combin.} \textbf{14} (2007), 23 pp.
	\bibitem{MSW13} G. Musiker, R. Schiffler, and L. Williams, Bases for cluster algebras from surfaces, \textsl{Compos. Math} \textbf{149} (2013), 2891--2944.
	\bibitem{N11} T. Nakanishi, Periodicities in cluster algebras and dilogarithm identities, in A. Skowronski \& K. Yamagata (eds.), Representations of algebras and related topics, EMS Series of Congress Reports, Eur. Math. Soc. (2011) 407--444.
	\bibitem{OST10} V. Ovsienko, R. Schwartz, and S. Tabachnikov, The pentagram map: a discrete integrable system, \textsl{Comm. Math. Phys.} \textbf{299} (2010), 409--446.
	\bibitem{OST13} V. Ovsienko, R. Schwartz, and S. Tabachnikov, Liouville-Arnold integrability of the pentagram map on closed polygons, \textsl{Duke Math. J.} \textbf{162} (2013), 2149--2196.
  	\bibitem{P15} P. Pylyavskyy, Zamolodchikov integrability via rings of invariants, arXiv preprint: 1506.05378v4.
	\bibitem{Q12} F. Qin, Quantum cluster variables via Serre polynomials, \textsl{J. Reine Angew. Math.} \textbf{668} (2012), 149--190.
  	\bibitem{R11} D. Rupel, On a quantum analogue of the Caldero-Chapoton formula, \textsl{Int. Math. Res. Not.} (2011), no. 14, 3207--3236.
  	\bibitem{R15} D. Rupel, Quantum cluster characters for valued quivers, \textsl{Trans. Amer. Math. Soc.} \textbf{367} (2015), 7061--7102.
	\bibitem{S92} R. Schwartz, The pentagram map, \textsl{Experiment. Math.} \textbf{1} (1992), 71--81.
	\bibitem{S08} R. Schwartz, Discrete monodromy, pentagrams, and the method of condensation, \textsl{J. Fixed Point Theory Appl.} \textbf{3} (2008), 379--409.
	\bibitem{S13} F. Soloviev, Integrability of the Pentagram Map, \textsl{Duke Math. J.} \textbf{162} (2013), 2815--2853.
	\bibitem{V07} A. Volkov, On the periodicity conjecture for $Y$-systems, \textsl{Comm. Math. Phys.} \textbf{276} (2007), 509--517.
  	\bibitem{W13} L. Williams, Cluster algebras: an introduction, arXiv preprint: 1212.6263v3.
\end{thebibliography}
\end{document}